\documentclass[12pt]{amsart}
\usepackage{amsmath,amsfonts,textcomp,longtable,pstricks,amscd,mathrsfs}%xypic,
\usepackage[dvips]{graphicx}

\newtheorem*{thm*}{Theorem}
\newtheorem{thm}[subsection]{Theorem }

\newtheorem{lemma}[subsection]{Lemma }
\newtheorem*{lemma*}{Lemma}

\newtheorem{sublemma}[subsubsection]{Lemma }

\newtheorem*{prop*}{Proposition}
\newtheorem{prop}[subsection]{Proposition }

\newtheorem*{corollary*}{Corollary}
\newtheorem{corollary}[subsection]{Corollary }

\theoremstyle{definition}

\newtheorem*{rem*}{Remark}

\newtheorem*{ex*}{Example}
\newtheorem{ex}[subsection]{Example }

\def\PP{{\mathbb P}}

\def\QQ{{\mathbb Q}}
\def\ZZ{{\mathbb Z}}
\def\NN{{\mathbb N}}

 \DeclareMathOperator{\Hom}{Hom}

\DeclareMathOperator{\D}{D}

\DeclareMathOperator{\Tot}{Tot}

 \DeclareMathOperator{\tr}{tr}
 
 \DeclareMathOperator{\Spec}{Spec}
 
\DeclareMathOperator{\Tr}{Tr}

\DeclareMathOperator{\cone}{cone}
\DeclareMathOperator{\Gr}{Gr}

\def\OO{\mathcal{O}}
\def\UU{\mathcal{U}}

\renewcommand \phi {\varphi}
\renewcommand \rho {\varrho}

\DeclareMathAlphabet{\mathpzc}{OT1}{pzc}{m}{it}

\newcommand{\DG}{\mathrm{DG}}
\newcommand{\Hot}{\mathrm{Hot}}
\newcommand{\Acycl}{\mathrm{Acycl}}

\newcommand{\cHom}{{\mathscr H}\mkern -4mu\mathit{om}}
\newcommand{\RHom}{R{\mathscr{H}\mkern -4mu\mathit{om}}}
\newcommand{\cExt}{{\mathscr E}\mkern -3mu\mathit{xt}}
\newcommand{\cTor}{{\mathscr T}\mkern -4mu\mathit{or}}
\newcommand{\cDR}{{\mathscr{DR}}}
\newcommand{\DR}{\mathrm{DR}}
\newcommand{\DD}{\mathscr D}

\newcommand{\sinj}{{\mathrm{inj}}}

\newcommand{\sfl}{{\mathrm{fl}}}
\newcommand{\sperf}{{\mathrm{perf}}}
\newcommand{\sthick}{{\mathrm{thick}}}

\newcommand{\co}{{\mathrm{co}}}

\newcommand{\abs}{{\mathrm{abs}}}
\newcommand{\gr}{{\mathrm{gr}}}

\newcommand{\sqc}{{\mathrm{qcoh}}}
\newcommand{\scoh}{{\mathrm{coh}}}
\newcommand{\st}{{\mathrm{st}}}

\def\qc#1{\operatorname{\mathrm{\Omega_{#1}--\,\sqc}}}

\def\Dmod#1{\operatorname{\mathrm{\DD_{#1}--\,\sqc}}}

\def\coh#1{\operatorname{\mathrm{\Omega_{#1}--\,\scoh}}}
\def\Dcoh#1{\operatorname{\mathrm{\DD_{#1}--\,\scoh}}}
\def\inj#1{\operatorname{\mathrm{\Omega_{#1}--\,\sqc_\sinj}}}
\def\fl#1{\operatorname{\mathrm{\Omega_{#1}--\,\sqc_\sfl}}}

\def\perf#1{\operatorname{\mathrm{\Omega_{#1}--\,\sperf}}}

\def\thick#1{\operatorname{\mathrm{\Omega_{#1}--\,\sqc_\sthick}}}
\def\d{\partial}

\begin{document}
\author{Sergey Rybakov}
\thanks{The author is partially supported by AG Laboratory GU-HSE, RF government grant, ag. 11 11.G34.31.0023, 
and by RFBR grants no. 11-01-12072 and 11-01-00395}
\address{Poncelet laboratory (UMI 2615 of CNRS and Independent University of
Moscow)}
\address{Institute for information transmission problems of the Russian Academy of Sciences}
\address{Laboratory of Algebraic Geometry, NRU HSE, 7 Vavilova Str., Moscow, Russia, 117312}
\email{rybakov@mccme.ru, rybakov.sergey@gmail.com}%

\title{DG-modules over de Rham DG-algebra}
\date{}


\begin{abstract}
For a morphism of smooth schemes over a regular affine base we define functors of derived direct image and extraordinary inverse image on coderived categories of DG-modules over de Rham DG-algebras. Positselski proved that for a smooth algebraic variety $X$ over a field $k$ of characteristic zero the coderived category of DG-modules over $\Omega^\bullet_{X/k}$ is equivalent to the unbounded derived category of quasi-coherent right $\DD_X$-modules. We prove that our functors correspond to the functors of the same name for $\DD_X$-modules under Positselski equivalence. 
\end{abstract}
\maketitle

\section{Introduction}
Let $X$ be a smooth algebraic variety over a field $k$ of characteristic zero. Its de Rham complex $\Omega_X=\Omega^\bullet_{X/k}$ is a sheaf of DG-algebras, and quasi-coherent sheaves of DG-modules over $\Omega_X$ (we call them $\Omega_X$-modules) form a DG-category, but the derived category of this category is usually not interesting. Recently Positselski defined the coderived category $\D^\co(\qc{X})$ of quasi-coherent $\Omega_X$-modules and proved that the unbounded derived category of quasi-coherent left $\DD_X$-modules $\D(\Dmod{X})$ and the coderived category $\D^\co(\qc{X})$ are equivalent~\cite{P2}. This result is a generalization of Koszul duality for symmetric and exterior algebras. In this paper we construct functors of derived tensor product, direct image and extraordinary inverse image on $\Omega_X$-modules and prove that these functors correspond to the functors of the same name for $\DD_X$-modules under Positselski equivalence (see section~\ref{comparison}). For coherent $\Omega_X$-modules we define the duality functor which correspond to duality on $\DD_X$-modules. Before this, we prove some preliminary results which have obvious $\DD_X$-module counterparts.

Some of these constructions are well-known for DG-categories of $\Omega$-modules (for example,~\cite{Ka} and~\cite[7.2]{BD}), but as far as author knows there was no attempt to write down this definitions for coderived categories of $\Omega$-modules. We feel it is important to do it.

%The paper is organized as follows. In section $2$ we establish notation and recall some basic facts on the DG-category of $\Omega_X$-modules. In section $3$ we collect results on exotic derived categories due to Positselski. %Of course, any error in this section is an author's fault. 
%Section $4$ contains definitions and basic results on the functors we are interested in: derived tensor product, internal $\Hom$, direct image and extraordinary inverse image. 
%Section~\ref{algebra} contains some technical results on injective hulls and on the stupid filtration on  $\Omega_X$-modules.
%In section~\ref{kt} we prove an $\Omega_X$-module version of Kashivara theorem, and in section~\ref{tce} we construct the trace isomorphism. %The last construction involves an exotic type of Koszul resolution which seems to be new.
%Section~\ref{sec_coh} is devoted to coherent $\DD_X$-modules. 
%Finally, in section~\ref{comparison} we prove the main comparison result: Positselski equivalence takes derived tensor product, direct image and extraordinary inverse image on $\Omega_X$-modules, and duality on coherent $\Omega_X$-modules to corresponding functors on quasi-coherent $\DD_X$-modules.

I am very grateful to L. Positselski for explaining to me his results on exotic derived categories and for many useful conversations. I also thank Sergey Arkhipov and Sergey Gorchinskiy for helpful discussions on the topic. 

\section{Quasi-coherent sheaves over $\Omega_X$}
\subsection{DG-rings and DG-modules.}
A \emph{DG-ring} is a graded ring $R=\oplus_{i\in\ZZ}R^i$ endowed with an odd derivation $d: R\to R$ of degree $1$ such that $d^2=0$. We write $R^\gr$ when we view $R$ as a graded ring (or a DG-ring with zero diferential). We denote by $|x|$ the degree of a homogeneous element $x\in R^{|x|}$.
A DG-ring is called \emph{commutative} if $xy=(-1)^{|x||y|}yx$ for any pair of homogeneous elements $x$ and $y$ of $R$. In what follows we assume that $R$ is a commutative DG-ring.

A (left) \emph{DG-module} $M$ over $R$ is a graded $R$-module $M=\oplus_{i\in\ZZ}M^i$ endowed with differential $d_M:M\to M$ of degree $1$ such that $d_M^2=0$, and for any pair of homogeneous elements $x\in R$ and $m\in M$ the Leibniz rule holds: $$d_M(xm)=dx\cdot m+(-1)^{|x|}xd_M(m).$$ The module $M[j]=\oplus_{i\in\ZZ}M^{i+j}$ is called the \emph{shift} of $M$ by an integer $j$. The differential on $M[j]$ is defined as $(-1)^jd_M$.

Let $M$ and $N$ be DG-modules over $R$. The \emph{DG-module of homomorphisms} $\Hom_R(M,N)$ from $M$ to $N$ over $R$ is defined as follows. The component $\Hom_R^i(M,N)$ consists of all homogeneous maps $f:M\to N$ of degree $i$ such that $f(xm)=(-1)^{i|x|}xf(m)$ for homogeneous $x\in R$ and $m\in M$. Since $R$ is commutative, the formula $(xf)(m)=xf(m)$ defines an action of $R$ on $\Hom_R(M,N)$. The differential in the complex $\Hom_R(M,N)$ is given by the formula: $df(m)=d_N(f(m))-(-1)^{|f|}f(d_M(m))$. Clearly, $d^2(f)=0$, and for any composable morphisms $f$ and $g$ one has $d(fg)=df\cdot g+(-1)^{|f|}fdg$. %Clearly, $\Hom_R(M,N)$ is a DG-module over $R$. 

We define \emph{tensor product DG-module} $M\otimes_R N$ of DG-modules $M$ and $N$ over $R$ as the tensor product $M^\gr\otimes_{R^\gr} N^\gr$ with the differential given by the formula: $d(m\otimes n)=dm\otimes n+(-1)^{|m|}m\otimes dn$.

\subsection{DG-modules over the de Rham DG-algebra.}
In this paper we always assume that all schemes are smooth, quasi-projective and connected over a regular connected affine base $S$. We denote by $\dim X$ the Krull dimension of $X$, and by $\Omega_X$ the sheaf of de Rham DG-algebras $\Omega^\bullet_{X/S}$. Note that $\Omega_X$ is a sheaf of commutative DG-algebras over $\OO_X$. In what follows sheaves of DG-modules over $\Omega_X$ are assumed to be quasi-coherent over $\OO_X$. We call them $\Omega_X$-modules for short. $\Omega_X$-modules form a DG-category $\qc{X}$ with shifts, cones and arbitrary direct sums.
The \emph{complex of homomorphisms} $\Hom_{\Omega_X}(M,N)$ from $M$ to $N$ over $\Omega_X$ is defined in a similar way as in the previous section. The component $\Hom_{\Omega_X}^i(M,N)$ consists of all homogeneous maps $f:M\to N$ of degree $i$ such that over any affine open $U\subset X$ and for any $\omega\in\Omega_U$ and $m\in M(U)$ the following relation holds: $f(\omega m)=(-1)^{i|\omega|}\omega f(m)$. The formula for the differential is the same as before. Note that $\Hom_{\Omega_X}^\bullet(M,N)$ is a complex of $\OO_S(S)$-modules. 

An $\Omega_X$-module $M$ is called \emph{coherent} if $M$ is coherent as an $\OO_X$-module. In particular, there is a finite set of indexes $i$ such that the graded component $M^i$ is not zero. Coherent $\Omega_X$-modules form a DG-category $\coh{X}$ with shifts, cones and finite direct sums.

Let $M$ be a coherent $\Omega_X$-module. In this case we can define an \emph{$\Omega_X$-module $\cHom_{\Omega_X}(M,N)$ of homomorphisms} from $M$ to $N$ over $\Omega_X$. For any open $U\subset X$ the DG-module $\cHom_{\Omega_X}(M,N)(U)$ is the DG-module of homomorphisms $\Hom_{\Omega_U}(M(U),N(U))$. Since $M$ is coherent, $\cHom_{\Omega_X}(M,N)$ is quasi-coherent.

Let $M$ and $N$ be $\Omega_X$-modules. One defines the \emph{tensor product} $\Omega_X$-module $M\otimes_{\Omega_X} N$ as the associated sheaf to the presheaf $U\mapsto M(U)\otimes_{\Omega_U} N(U)$. 
%The module $(M\otimes_{\Omega_X} N)(U)$ can be obtained from DG-modules over $\Omega_U$ by taking limit of an injective system of DG-modules and taking kernel of a morphism of DG-modules. Thus $(M\otimes_{\Omega_X} N)(U)$ is a DG-module over $\Omega_U$, and $(M\otimes_{\Omega_X} N)$ is an $\Omega_X$-module.

We define a \emph{closed morphism} $f$ as an element of $\Hom_{\Omega_X}^0(M,N)$ such that $d(f)=0$. 
%A morphism $f:M\to N$ is called \emph{nill-homotopic}, if $f=d(h)$ for some $h\in\Hom_{\Omega_X}^{-1}(M,N)$.
The $\Omega_X$-module $M$ is called \emph{contractible} if the identity morphism is homotopic to zero, i.e. $Id_M=d(h)$ for some $h\in\Hom_{\Omega_X}^{-1}(M,N)$ .
Let $f:M\to N$ be a closed morphism of $\Omega_X$-modules. The \emph{cone} of $f$ is the following $\Omega_X$-module. As a graded module $\cone(f)=M[1]\oplus N$, and the differential is given by the formula:$$d_{\cone(f)}(m,n)=(-d_M(m),f(m)+d_N(n)).$$
There is an exact sequence of $\Omega_X$-modules with closed differentials:
$$0\to N\to\cone(f)\to M[1]\to 0.$$

We need the following generalization of the cone. Let $\dots\to M_i\to M_{i-1}\to\dots$ be a complex of $\Omega_X$-modules with closed differentials $d_i:M_i\to M_{i-1}$. Take the graded module $\Tot^{\oplus}(M_\bullet)=\oplus_i M_i[i]$ with differential $$d(m)=(-1)^id_{M_i}(m)+d_i(m)$$ for $m\in M_i$. The $\Omega_X$-module $\Tot^{\oplus}(M_\bullet)$ is called the \emph{total object} of the complex $M_\bullet$ formed by (infinite) direct sums.

Let $D$ be a DG-category with shifts and cones. 
The \emph{homotopy category} $\Hot(D)$ is the additive category with the same class of objects as $D$ and groups of morphisms given by $\Hom_{\Hot(D)}(X,Y)=H^0(\Hom_D(X,Y))$.
The homotopy category $\Hot(D)$ is a triangulated category~\cite{BK}.

\section{Positselski results on coderived categories}
This section is an overview of some results on exotic derived categories due to Positselski. We present here a down to earth treatment of the subject~\cite{PLJ}. The reader can find the results in full generality in the papers~\cite{P1,P2,P3,PP}.

\subsection{} Fix an exact subcategory $E$ in the abelian category of quasi-coherent graded $\Omega^\gr_X$-modules. We say that $\Omega_X$-module is an \emph{$E$-module} if its  underlying graded module belongs to $E$. An $\Omega_X$-module is called \emph{absolutely acyclic} with respect to $E$ if it belongs to the minimal thick subcategory of the homotopy category $\Hot(\qc{X})$ containing total objects of exact triples of $E$-modules. We denote this subcategory by $\Acycl^\abs_E$. The quotient category of the homotopy category of $E$-modules by $\Acycl^\abs_E$ is called \emph{the absolute derived category of $E$-modules}. 

For example, let $E$ be the subcategory of coherent $\Omega_X^\gr$-modules. Then $\D^\abs(\coh{X})=\Hot(\coh{X})/\Acycl^\abs_\scoh$ is the absolute derived category of
coherent $\Omega_X$-modules. 

Suppose that $E$ is closed under arbitrary direct sums. An $\Omega_X$-module is called \emph{coacyclic} with respect to $E$ if it belongs to the minimal thick subcategory of the homotopy category $\Hot(\qc{X})$ containing total objects of exact triples of $E$-modules and closed under arbitrary direct sums. We denote this subcategory by $\Acycl^\co_E$. The quotient category of the homotopy category of $E$-modules by $\Acycl^\co_E$ is called \emph{the coderived category} of $E$-modules. 

If $E$ is the category of quasi-coherent $\Omega_X^\gr$-modules, we call coacyclic $\Omega_X$-modules with respect to $E$ \emph{coacyclic}, and denote the subcategory of coacyclic $\Omega_X$-modules by $\Acycl^\co$.
Similarly, the category of \emph{absolutely acyclic} $\Omega_X$-modules $\Acycl^\abs$ is the category of absolutely acyclic $\Omega_X$-modules with respect to $E$.
The category $\D^\co(\qc{X})=\Hot(\qc{X})/\Acycl^\co$ is called the coderived category of $\Omega_X$-modules. We call two $\Omega_X$-modules \emph{co-isomorphic} if they are isomorphic in $\D^\co(\qc{X})$. A curious example of co-isomorphism is given in Lemma~\ref{coiso}.
%co-isomorphism of $\Omega_Y$-modules $\OO_Y\to\Omega_Y\otimes_{\OO_Y}\omega_Y^{-1}\otimes_{\OO_Y}\DD_Y[\dim Y]$

A coherent absolutely acyclic $\Omega_X$-module is absolutely acyclic with respect to the category of coherent $\Omega_X^\gr$-modules, and the natural functor from $\D^\abs(\coh{X})$ to $\D^\co(\qc{X})$ is full and faithful~\cite[3.11]{P2}.

\begin{lemma}\label{colemma}
\begin{enumerate}
\item Let $M=\varinjlim_{n\in \NN}M_n$ be the limit of an inductive system of coacyclic $\Omega_X$-modules with closed morphisms. Then $M$ is coacyclic. 
\item Let $M_0\xrightarrow{d_0}M_1\to \dots$ be an exact sequence of $\Omega_X$-modules with closed differentials. Then $M=\ker d_0$ and $\Tot^\oplus M_\bullet$ are co-isomorphic. 
\item Let $M=\cup_{n\in\NN} M_n$ be an $\Omega_X$-module with an increasing filtration such that $M_{n+1}/M_n$ is coacyclic for any $n\in\NN$. Then $M/M_0$ is coacyclic.
\end{enumerate}
\end{lemma}
\proof The first assertion follows from definitions since $M=\varinjlim_{n\in \NN}M_n$ is the cone of the injective map $id-d:\oplus_{n\in\NN}M_n\to\oplus_{n\in\NN}M_n$. To prove the second assertion note that the total object of the sequence $0\to M\to M_0\to\dots\to M_n\to\ker d_{n+1}\to 0$ is absolutely acyclic. Hence the direct limit $N$ of such total objects is coacyclic, and we have an exact sequence of $\Omega_X$-modules:
$$0\to \Tot^\oplus M_\bullet\to N\to M[-1]\to 0.$$ The third assertion follows from the fact that $M/M_0=\varinjlim_{n\in \NN}M_n/M_0$.
\qed

\subsection{}\label{dg} An $\Omega_X$-module $I$ is called \emph{injective} if $I$ is injective as a graded $\Omega_X^\gr$-module. By~\cite[Theorem A3]{P3}, $I$ is injective if and only if the restriction of $I$ to any open subscheme $U$ is an injective $\Omega_U$-module.
An $\Omega_X$-module $F$ is called \emph{flat} if $F$ is flat as a graded $\Omega_X^\gr$-module. Let $M$ be an $\Omega_X$-module. We would like to find an injective $\Omega_X$-module $I$ and a flat $\Omega_X$-module $F$ which are co-isomorphic to $M$. As a first approximation to $I$ and $F$ we construct an injective $\Omega_X$-module $I_0$ with an injection $M\to I_0$ and a locally free (over $\Omega_X^\gr$) $\Omega_X$-module $P_0$ and a surjection $P_0\to M$.

There is an injective morphism of graded $\OO_X$-modules $M\to J$ to an injective quasi-coherent graded $\OO_X$-module $J$ and a surjective $\OO_X$-morphism $Q\to M$ from a locally free graded $\OO_X$-module $Q$. We are going to construct functors $\DG^+$ and $\DG^-$ from graded $\OO_X$-modules to $\Omega_X$-modules such that $I_0=\DG^-(J)$, and $P_0=\DG^+(Q)$.

\subsubsection{}\label{dg-} Let $J$ be a quasi-coherent $\OO_X$-module. The module $\DG^-(J)$ is cofreely generated by $J$ over $\Omega_X$. Namely, take a coinduced module $J'=\cHom_{\OO_X}(\Omega_X^\gr,J)$. It is a graded $\Omega_X^\gr$-module such that for any $\Omega_X^\gr$-module $N$ there is a natural isomorphism:
$$\Phi:\cHom_{\OO_X}(N,J)\to\cHom_{\Omega_X^\gr}(N,J').$$
The graded sheaf of abelian groups $\DG^-(J)$ is defined as follows. For an open $U\subset X$ the group $\DG^-(J)(U)$ consists of all formal expressions of the form $d^{-1}x+y$, where $x$, $y\in J'(U)$, and $\deg d^{-1}x=\deg x-1$. The differential and the action of $\Omega_U$ on $\DG^-(J)(U)$ is given by the formulas: $d(d^{-1}x+y)=x$ and $\omega(d^{-1}x+y)=\omega y+d^{-1}(d(\omega) y)+(-1)^{|\omega|}d^{-1}(\omega x)$ for $\omega\in\Omega_U$.

For any $\Omega_X$-module $M$ there is a bijective correspondence between morphisms of graded $\OO_X$-modules $f:M\to J$ and closed morphisms of $\Omega_X$-modules $g:M\to\DG^-(J)$ which is described by the formula: $g(m)=d^{-1}(\Phi(f)(dm))+\Phi(f)(m)$. The module $\DG^-(J)$ is contractible. If $J$ is an injective $\OO_X$-module, then $\DG^-(J)\cong J'\oplus J'[-1]$ as a graded $\Omega_X^\gr$-module. 

\subsubsection{}\label{dg+} The module $\DG^+(Q)$ is constructed in a similar way.
As a graded sheaf of abelian groups $$\DG^+(Q)(U)=\{x+dy:x,y\in (Q\otimes_{\OO_X}\Omega_X)(U)\},$$ where $\deg dy=\deg y+1$. The differential and the action of $\Omega_U$ on $\DG^+(Q)(U)$ is given by the formulas: $d(x+dy)=dx$ and $\omega(x+dy)=\omega x-(-1)^{|\omega|}(d\omega y)+(-1)^{|\omega|}d(\omega y)$ for $\omega\in\Omega_U$.

As before, there is a bijective correspondence between morphisms of graded $\OO_X$-modules $f:Q\to M$ and closed morphisms of $\Omega_X$-modules $g:\DG^+(Q)\to M$, and $\DG^+(Q)$ is contractible.
If $Q$ is locally free over $\OO_X$, then $\DG^+(Q)$ is locally free as a
%$\cong Q\otimes_{\OO_X}\Omega_X^\gr\oplus Q\otimes_{\OO_X}\Omega_X^\gr[1]$ as a 
graded $\Omega_X^\gr$-module.

\subsection{Injective resolutions.}\label{inj_resolution}
Denote the category of injective $\Omega_X$-modules by $\inj{X}$.

\begin{thm*}
The coderived category $\D^\co(\qc{X})$ is equivalent to the homotopy category $\Hot(\inj{X})$.
\end{thm*}
\proof Let $M$ be an $\Omega_X$-module. We construct an {\it injective resolution} $I$ of $M$ as follows.
By section~\ref{dg}, we have an inclusion $M\to I_0$, where $I_0$ is injective. A common argument gives an exact complex $I_\bullet=(I_0\to I_1\to\dots)$ of injective $\Omega_X$-modules with closed differentials. By Lemma~\ref{colemma}.1, $I=\Tot^\oplus I_\bullet$ is co-isomorphic to $M$. Since $\Omega_X^\gr$ is Noetherian, the direct sum of injective modules is injective, thus $I$ is injective.

By \cite[Theorem 3.5]{P2}, for any coacyclic $\Omega_X$-module $L$ the complex $\Hom_{\Omega_X}(L,I)$ is acyclic. Thus $\Hot(\inj{X})$ and coacyclic $\Omega_X$-modules form a semiorthogonal decomposition of $\Hot(\qc{X})$. By \cite[Lemma 1.3]{P2}, this proves the theorem.\qed

\subsection{Thick $\Omega_X$-modules.}\label{thick_mod}
We call an $\Omega_X$-module $A$ \emph{(co)induced} if $A$ is (co)induced as an $\Omega_X^\gr$-module, i.e., there exists an $\OO_X$-module $L$ such that $A^\gr\cong L\otimes_{\OO_X}\Omega_X^\gr$ (or $A^\gr\cong\cHom_{\OO_X}(\Omega_X^\gr,L)$). Note that any induced $\Omega_X^\gr$-module is coinduced and vice versa.

\begin{thm*} Let $A$ be a quasi-coherent $\Omega_X^\gr$-module.
The following conditions are equivalent:
\begin{enumerate}
\item $\cTor^i_{\Omega_X^\gr}(\OO_X,A)=0$ for any $i>0$;
\item $\cTor^i_{\Omega_X^\gr}(F,A)=0$ for any $i>0$, and for any $\OO_X$-flat $\Omega_X^\gr$-module $F$;
\item $\cExt^i_{\Omega_X^\gr}(A,J)=0$ for any $i>0$, and for any $\OO_X$-injective $\Omega_X^\gr$-module $J$ such that $\Omega^1_XJ=0$;
\item $\cExt^i_{\Omega_X^\gr}(A,I)=0$ for any $i>0$, and for any $\OO_X$-injective $\Omega_X^\gr$-module $I$;
\item $\cExt^i_{\Omega_X^\gr}(\OO_X,A)=0$ for any $i>0$;
\item $\cExt^i_{\Omega_X^\gr}(P,A)=0$ for any $i>0$, and for any %perfect 
$\Omega_X^\gr$-module $P$ which is locally projective over $\OO_X$.
\end{enumerate}
If $A$ is (co)induced, then $(1)-(6)$ holds.
\end{thm*}
\proof
Clearly, $(2)\Rightarrow(1)$, $(4)\Rightarrow(3)$, $(6)\Rightarrow(5)$, and if $A$ is (co)induced, then $(1)-(6)$ holds. 

Let us prove that $(1)\Leftrightarrow(3)$. Take a locally free over $\Omega_X^\gr$ resolution for $A$:
$$\dots\to P_n\to\dots\to P_0\to A\to 0.\eqno (*)$$
The complex $$\dots\to P_n\otimes_{\Omega_X^\gr}\OO_X\to\dots\to P_0\otimes_{\Omega_X^\gr}\OO_X\to 0$$ 
computes $\cTor^i_{\Omega_X^\gr}(\OO_X,A)$.
%is a resolution for $A\otimes_{\Omega_X^\gr}\OO_X$ over $\OO_X$. 
Apply $\cHom_{\OO_X}(\cdot,J)$ to this complex. We obtain a complex which computes $\cExt^\bullet_{\Omega_X^\gr}(A,J)$. Since $J$ is an arbitrary $\OO_X$-injective module, the groups $\cExt^i_{\Omega_X^\gr}(A,J)$ are trivial if and only if $\cTor^i_{\Omega_X^\gr}(\OO_X,A)=0$.

$(4)\Rightarrow(2)$. Let $J$ be an injective $\OO_X$-module. The module $I=\cHom_{\OO_X}(F,J)$ is $\OO_X$-injective $\Omega_X^\gr$-module. Apply $\cHom_{\Omega_X^\gr}(\cdot,I)$ to the resolution $(*)$. If $(4)$ holds, we get a resolution for $\cHom_{\OO_X}(A\otimes_{\Omega_X^\gr}F,J)$. Since $J$ is an arbitrary injective $\OO_X$-module, $\cTor^i_{\Omega_X^\gr}(F,A)=0$ for any $i>0$.

$(1)\Rightarrow(6)$. First, we prove a lemma.
\begin{lemma*}
Assume $(1)$. Then $A$ has a finite left resolution over $\Omega_X^\gr$ which is locally (co)induced.
\end{lemma*}
\proof
Take the truncation of the resolution $(*)$:
$$0\to Q\to P_n\to\dots\to P_0\to A\to 0, \eqno (**)$$ where $n=\dim X$.
Then $Q$ is locally projective $\OO_X$-module. Tensor $(**)$ with $\OO_X$ over $\Omega_X^\gr$. By condition $(1)$, the module $L=Q\otimes_{\Omega_X^\gr}\OO_X$ is also locally projective over $\OO_X$, and $\cTor^i_{\Omega_X^\gr}(\OO_X,Q)=0$ for any $i>0$. We have the natural surjective morphism of $\Omega_X^\gr$-modules $$Q\to L\cong Q/\Omega^1_XQ.$$ This implies that locally there is a morphism of $\OO_X$-modules $L\to Q$ which induces a surjective morphism of $\Omega_X^\gr$-modules $$\phi:L\otimes_{\OO_X}\Omega_X^\gr\to Q.$$ Let $K=\ker\phi$. Clearly, $\cTor^i_{\Omega_X^\gr}(\OO_X,K)=0$ for any $i\geq 0$, and thus $K=0$. We proved that locally $Q$ is a (co)induced module.
\qed

By the Lemma, the module $A$ has a resolution $$0\to P_n\to\dots\to P_0\to A\to 0$$ such that for any $P_k$ the property $(6)$ holds: $\cExt^i_{\Omega_X^\gr}(P,P_k)=0$ for any $i>0$. It follows that $(6)$ is true for $A$. 
To prove $(5)\Rightarrow(4)$ one uses almost the same argument involving right resolution of $A$ instead of the resolution $(*)$. 
\qed

We call an $\Omega_X^\gr$-module $A$ \emph{thick} if conditions $(1)-(6)$ holds for $A$. We call $\Omega_X$-module \emph{thick} if it is thick as a graded $\Omega_X^\gr$-module. Positselski call such modules \emph{weakly relatively admissible}. The following result can be proved using the argument given in the previous Lemma.

\begin{corollary*} Let $A$ be a thick $\Omega_X^\gr$-module. Suppose that
$A$ is locally projective (flat or injective) over $\OO_X$. Then $A$ is locally projective (respectively flat or injective) over $\Omega_X^\gr$.
\end{corollary*}

\subsection{Thick resolutions.}\label{thick_resolution}
\begin{thm*}
The coderived category $\D^\co(\qc{X})$ is equivalent to the absolute derived category of $\D^\abs(\thick{X})$ of thick $\Omega_X$-modules.
\end{thm*}
\proof
We prove that subcategory of absolutely acyclic thick $\Omega_X$-modules and the homotopy category of injective $\Omega_X$-modules form a semiorthogonal decomposition of the homotopy category of thick $\Omega_X$-modules. The theorem then follows from \cite[Lemma 1.3]{P2} and Theorem~\ref{inj_resolution}.

Let $A$ be a thick $\Omega_X$-module. Find an $\OO_X$-injection of $A$ into an injective $\OO_X$-module $J_0$. By~\ref{dg} we obtain an inclusion $A\to I_0=\DG^-(J_0)$ into a coinduced injective $\Omega_X$-module. In this way we construct first terms of the resolution for $A$ by coinduced injective $\Omega_X$-modules:
$$I_0\to I_1\to\dots\to I_n,$$ where all morphisms are closed. Let $J$ be the cokernel of the map $I_{n-1}\to I_n$. Then $J$ is thick. Moreover,  if $n\geq\dim X$, then $J$ is injective over $\OO_X$. By Corollary~\ref{thick_mod}, $J$ is an injective $\Omega_X$-module. It follows that $$I=\Tot^\oplus(I_0\to\dots\to I_n\to J)$$ is injective too. Moreover, the cokernel $L$ of the natural morphism $A\to I$ is an absolutely acyclic thick $\Omega_X$-module. By \cite[Theorem 3.5]{P2}, $\Hom_{\Omega_X}(L,I)$ is acyclic (see also Lemma~\ref{hom_tensor}).
\qed

\begin{rem*}%\label{free_resolution}
In fact, we proved that any $\Omega_X$-module $M$ has a \emph{coinduced resolution}, i.e there exists a coinduced $\Omega_X$-module $A$ and a morphism $M\to A$ with absolutely acyclic cone. 
\end{rem*}

\subsection{Flat resolutions.}\label{flat_resolution}
\begin{thm*} The coderived category $\D^\co(\qc{X})$ is equivalent to the absolute derived category  $\D^\abs(\fl{X})$ of flat $\Omega_X$-modules.
\end{thm*}
\proof 
By Theorem~\ref{thick_resolution}, $\D^\co(\qc{X})$ is equivalent to $\D^\abs(\thick{X})$. Let us prove that the latter is equivalent to $\D^\abs(\fl{X})$. Let $A$ be a thick $\Omega_X$-module. As before, using~\ref{dg} find a finite exact complex $$0\to P_n\to\dots\to P_0\to A\to 0$$ of  $\Omega_X$-modules such that $P_0$, \dots, $P_{n-1}$ are flat (and even locally free over $\Omega_X^\gr$) and differentials are closed. If $n$ is greater than dimension of $X$, the module $P_n$ is also flat by Corollary~\ref{thick_mod}. It follows that $P=\Tot^{\oplus} P_{\bullet}$ is flat and $\cone(P\to A)$ is absolutely acyclic with respect to $\thick{X}$.

We proved that the natural functor $\D^\abs(\thick{X})\to \D^\abs(\fl{X})$ is essentially surjective. We have to show that it is full and faithful. The key point is to prove that any object of $\fl{X}$ that is absolutely acyclic with respect to $\thick{X}$ is absolutely acyclic with respect to $\fl{X}$. This is done in~\cite[Section~3.2]{PP}. \qed

\subsection{Perfect resolutions.}\label{perf_resolution}
We need a one more type of resolutions. Denote the category of $\OO_X$-locally free coherent $\Omega_X$-modules by $\perf{X}$. We call such $\Omega_X$-modules \emph{perfect}. This is an exact category and we can define the absolute derived category $\D^\abs(\perf{X})$.

\begin{thm*} The absolute derived category $\D^\abs(\coh{X})$ of coherent $\Omega_X$-modules is equivalent to the absolute derived category $\D^\abs(\perf{X})$ of perfect $\Omega_X$-modules.
\end{thm*}
\proof Let $M$ be a coherent $\Omega_X$-module.
By section~\ref{dg}, there is a finite exact complex $0\to P_n\to\dots\to P_0\to M\to 0$ of  $\Omega_X$-modules such that $P_0$, \dots, $P_{n-1}$ are perfect (and even locally free over $\Omega_X^\gr$) and differentials are closed. If $n\geq\dim X$, the module $P_n$ is also perfect. It follows that $\Tot^{\oplus} P_{\bullet}$ is perfect and co-isomorphic to $M$. As before, one can prove that any object of $\perf{X}$ that is absolutely acyclic with respect to $\coh{X}$ is absolutely acyclic with respect to $\perf{X}$ (see~\cite[Section~3.2]{PP}). 
\qed

\begin{rem*}%\label{free_resolution}
For any $\Omega_X$-module $M$ there is a flat resolution $F$ which is locally free over $\Omega_X^\gr$. Indeed, the same construction as before gives a resolution $P$ of $M$ which is locally free over $\OO_X$, and the coinduced resolution of $P$ is locally free over $\Omega_X^\gr$.
\end{rem*}

\section{Inverse and direct images}
Let $f:X\to Y$ be a morphism of smooth schemes over $S$. Denote by $f^\bullet$ and $f_\bullet$ the inverse and direct image of sheaves of abelian groups associated to $f$. We have a morphism of sheaves of DG-algebras on $X$: $f^\bullet\Omega_Y\to\Omega_X$. This gives a morphism of DG-ringed spaces $(X,\Omega_X)\to (Y,\Omega_Y)$, and allows one to define inverse and direct image functors on coderived categories of $\Omega$-modules. 
%But first we establish some preliminary results on derived $\Hom$ and tensor product. 

\subsection{Derived $\Hom$ and tensor product.}\label{hom_tensor}
First we define functors of derived $\Hom$ and (shifted) tensor product:
$$\RHom_{\Omega_X}:\D^\abs(\coh{X})\times\D^\co(\qc{X})\to\D^\co(\qc{X}),$$
$$\otimes^D_{\Omega_X}:\D^\co(\qc{X})\times\D^\co(\qc{X})\to\D^\co(\qc{X}).$$

\begin{lemma*}
\begin{enumerate}
\item For any absolutely acyclic coherent $\Omega_X$-module $M$ and injective $\Omega_X$-module $I$ the $\Omega_X$-module $\cHom_{\Omega_X}(M,I)$ is absolutely acyclic.
\item For any coherent $\Omega_X$-module $M$ and coacyclic injective $\Omega_X$-module $I$ the $\Omega_X$-module $\cHom_{\Omega_X}(M,I)$ is coacyclic.
\item For any coacyclic $\Omega_X$-module $M$ and flat $\Omega_X$-module $F$ the $\Omega_X$-module $M\otimes_{\Omega_X}F$ is coacyclic.
\item Let $F$ be a flat absolutely acyclic {\rm(}with respect to the additive category of flat $\Omega_X$-modules{\rm)} $\Omega_X$-module. Then for any $\Omega_X$-module $M$ the $\Omega_X$-module $M\otimes_{\Omega_X}F$ is coacyclic.
\end{enumerate}
\end{lemma*}
\proof
The proof of $(1)$ is almost the same as~\cite[Theorem 3.5]{P2}. 
%Note that since $M$ is coherent, t commutes with arbitrary direct sums. As usual, it
Indeed, the functor $\cHom_{\Omega_X}(-,I)$ respects shifts, cones and direct sums.  Thus it is enough to check $(1)$ when $M$ is the total object of an exact triple of $\Omega_X$-modules. Clearly,
$\cHom_{\Omega_X}(M,I)$ is the the total object of an exact triple of $\Omega_X$-modules.
Similarly, if $I$ is co-acyclic, then $\Hom_{\Omega_X}(I,I)$ is acyclic, and $I$ is contractible. It follows that $\cHom_{\Omega_X}(M,I)$ is coacyclic. This proves $(2)$.

The functor $\cdot\otimes_{\Omega_X} F$ respects shifts, cones and arbitrary direct sums. We have to check $(3)$ when $M$ is the total object of an exact triple of $\Omega_X$-modules and $(4)$ when $F$ is the total object of an exact triple of flat $\Omega_X$-modules. Both times, $M\otimes_{\Omega_X}F$ is the the total object of an exact triple of $\Omega_X$-modules.
\qed

\subsubsection{}
We are ready to define the functor $\RHom_{\Omega_X}$. Take the restriction of the functor $$\cHom_{\Omega_X}:\Hot(\coh{X})\times\Hot(\qc{X})\to\Hot(\qc{X})$$ to the subcategory $\Hot(\coh{X})\times\Hot(\inj{X})$.
By Lemma $(1)$, this restriction factors through the functor 
\begin{equation}\label{hom_def}
\D^\abs(\coh{X})\times\Hot(\inj{X})\to\D^\co(\qc{X}).	
\end{equation}
By Theorem~\ref{inj_resolution}, $\Hot(\inj{X})$ is equivalent to $\D^\co(\qc{X})$. Combining this equivalence with the functor $(\ref{hom_def})$ we get the functor $\RHom_{\Omega_X}$.

\subsubsection{}
Let us define the derived tensor product. Let $n$ be the relative dimension of $X$ over $S$. By Lemma $(3)$ and $(4)$, the shifted tensor product $$\otimes_{\Omega_X}[2n]:\Hot(\qc{X})\times\Hot(\fl{X})\to\Hot(\qc{X})$$ factors through $$\D^\co(\qc{X})\times\D^\abs(\fl{X})\to\D^\co(\qc{X}).$$ By Theorem~\ref{flat_resolution}, $\D^\abs(\fl{X})$ is equivalent to $\D^\co(\qc{X})$, and we are done.

The tensor product can be obtained in two equivalent ways:
$$\D^\co(\qc{X})\times\D^\co(\qc{X})\to\D^\abs(\fl{X})\times\D^\co(\qc{X}) \xrightarrow{\otimes_{\Omega_X}[2n]}\D^\co(\qc{X}),$$
where the first arrow comes from the equivalence~\ref{flat_resolution}, or 
$$\D^\co(\qc{X})\times\D^\co(\qc{X})\to\D^\co(\qc{X})\times\D^\abs(\fl{X}) \xrightarrow{\otimes_{\Omega_X}[2n]}\D^\co(\qc{X}).$$

These definitions of tensor product are equivalent since both functors factors through $\D^\abs(\fl{X})\times\D^\abs(\fl{X})$.

\begin{thm}\label{thick_hom} 
Let $P$ be a perfect $\Omega_X$-module, and $A$ be a thick $\Omega_X$-module. Then
\begin{enumerate}
\item $\RHom_{\Omega_X}(P,A)\cong\cHom_{\Omega_X}(P,A)$, and
\item $P\otimes^D_{\Omega_X}A\cong P\otimes_{\Omega_X}A[2n]$
\end{enumerate}
in $\D^\co(\qc{X})$. 
\end{thm}
\proof
We prove $(1)$ and left $(2)$ to the reader. Take an injective resolution of $A$ from the proof of Theorem~\ref{thick_resolution}:$$I=\Tot^\oplus(I_0\to\dots\to I_n\to J).$$ We know that exact sequence $I_\bullet=(I_0\to\dots\to I_n\to J)$ gives an injective resolution of $A$ over $\Omega_X^\gr$. Since $P$ is perfect, the  complex of $\Omega_X^\gr$-modules $\cHom_{\Omega_X^\gr}(P^\gr,I^\gr_\bullet)$ has homology at zero term only. It follows that $\cHom_{\Omega_X}(P,I)$ is co-isomorphic to $\cHom_{\Omega_X}(P,A)$.
\qed

%\begin{rem*}
%We proved that for any thick module $A$ there are two functors: 
%$$\cHom_{\Omega_X}(\cdot,A):\D^\abs(\perf{X})\to\D^\co(\qc{X}),$$
%and $$\cdot\otimes_{\Omega_X}A[2n]:\D^\abs(\perf{X})\to\D^\co(\qc{X}).$$
%\end{rem*}

\subsection{Direct image.}\label{direct}
In this section we define the functor of direct image:
$$f_*:\D^\co(\qc{X})\to\D^\co(\qc{Y}).$$   
For any injective $\Omega_X$ module $I$ the direct image $f_\bullet(I)$ is an $\Omega_Y$-module since $f$ is a morphism of DG-ringed spaces $(X,\Omega_X)\to (Y,\Omega_Y)$. This gives the functor $$f_\circ:\Hot(\inj{X})\to\D^\co(\qc{Y}).$$ Any $\Omega_X$-module has an injective resolution, and we define $f_*$ as a composition of $f_\circ$ with the equivalence~\ref{inj_resolution}. It is well known that for an injective sheaf $I$ the sheaf $f_\bullet(I)$ is flabby, and $f_\bullet$ is exact on flabby sheaves. Thus, for a pair of composable maps $f$ and $g$ we have $(f\circ g)_*\cong f_*\circ g_*$.

By definition, $\Omega_S\cong\OO_S$, and $\Hom_{\Omega_X}(M,N)$ is a complex of $\OO_S(S)$-modules. Since $\OO_S$ has finite homological dimension and affine, $$\D^\co(\qc{S})\cong\D(\OO_S(S)\mathrm{-mod})$$ is an unbounded derived category of $\OO_S(S)$-modules~\cite[Remark 2.1]{P1}. %This fact allows one to compute $\Hom_{\D^\co(\qc{X})}(M,N)$ in terms of the DG-module $f_*\RHom_{\Omega_X}(M,N))$. 
%Since $S$ is affine, we easily deduce the following result.

\begin{prop*}
Let $f$ be the morphism from $X$ to $S$. And let $N\in\D^\co(\qc{X})$, and $M\in\D^\abs(\coh{X})$. Then
there is a canonical isomorphism in $\D(\OO_S(S)\mathrm{-mod})$: $$\Gamma(S,f_*\RHom_{\Omega_X}(M,N))\to \Hom_{\Omega_X}(M,I),$$
where $I$ is an injective resolution of $N$. In particular, the space of homomorphisms $\Hom_{\D^\co(\qc{X})}(M,N)$ is isomorphic to $H^0(S,f_*\RHom_{\Omega_X}(M,N))$.
%\Hom_{\Omega_X}(M,I))$.
\end{prop*}
\proof
The $\Omega_X$-module $\cHom_{\Omega_X}(M,I)$ is flabby. Thus $$\Gamma(S,f_*\RHom_{\Omega_X}(M,N))\cong\Gamma(S,f_\bullet\cHom_{\Omega_X}(M,I))\cong \Hom_{\Omega_X}(M,I).$$\qed

\subsection{Extraordinary inverse image.}\label{inverse}
Now, we define two inverse image functors: 
$$f^+,f^!:\D^\co(\qc{Y})\to\D^\co(\qc{X}).$$

By Theorem~\ref{flat_resolution}, it is enough to define the functors on $\D^\co(\fl{Y})$. We do it as follows:
$$f^+(F)=f^\bullet(F) \otimes_{f^\bullet\Omega_Y}\Omega_X,$$ where $F$ is flat.
{\it The extraordinary inverse image} is defined as $f^!(M)=f^+(M)[2(\dim X-\dim Y)]$.  
Clearly, for a pair of composable maps $f$ and $g$ we have $(f\circ g)^!\cong g^!\circ f^!$.
 
\begin{lemma*}
Suppose that $f:X\to Y$ is smooth. Then $\Omega_X$ is flat over $f^\bullet\Omega_Y$.
\end{lemma*}
\proof
It follows from the smoothness of $f$ that $\Omega_X^\gr$ is locally free over $f^\bullet(\Omega_Y^\gr)\otimes_{f^\bullet\OO_Y}\OO_X$,
and $\OO_X$ is flat over $f^\bullet\OO_Y$. Thus $\Omega_X^\gr$ is flat over  $f^\bullet(\Omega_Y^\gr)$.
\qed

\begin{corollary*} 
Suppose that $f:X\to Y$ is smooth, and $M$ is an $\Omega_Y$-module. Then
\begin{enumerate}
	\item $f^+(M)\cong f^\bullet(M)\otimes_{f^\bullet\Omega_Y}\Omega_X$;
	\item the extraordinary inverse image of a coherent $\Omega_Y$-module is coherent.
	\item If $f$ is an open embedding, then $f^!(M)\cong f^+(M)\cong f^\bullet(M)\otimes_{f^\bullet\OO_Y}\OO_X$.
\end{enumerate}
\end{corollary*}

\begin{thm}\label{adj_smooth}
Suppose that $f$ is smooth. Then functors $f_*$ and $f^+$ are adjoint. Let $N\in\D^\co(\qc{X})$, and let $M\in\D^\abs(\coh{X})$. Then
\begin{equation}\label{eq_cong_smooth}
\RHom_{\Omega_Y}(M,f_*N)\cong f_*\RHom_{\Omega_X}(f^+M,N).
\end{equation}
\end{thm}
\proof We prove the second assertion. The same argument proves that the functors are adjoint. 
Take an injective resolution $I$ for $N$. Since the functor $f^\bullet$ is exact, and functors $f_\bullet$ and $f^\bullet$ are adjoint, the module $f_\bullet(I)$ is injective.
 Moreover, homomorphism sheaf to an injective sheaf is flabby. Thus we have to prove that $$\cHom_{\Omega_Y}(M,f_\bullet I)\cong f_\bullet\cHom_{\Omega_X}(f^\bullet M\otimes_{f^\bullet\Omega_Y}\Omega_X,I).$$
And this can be shown in the usual way.
\qed

\begin{thm}{\rm(Projection formula.)}
Let $f:X\to Y$ be a morphism of smooth schemes over $S$. Let $M$ be an $\Omega_X$-module, and let $N$ be a $\Omega_Y$-module. Then $N\otimes^D_{\Omega_X}f_*M\cong f_*(f^!N\otimes^D_{\Omega_Y}M)$ in $\D^\co(\qc{Y})$.
\end{thm}
\proof
Let $I$ be an injective resolution of $M$, and $F$ be a flat resolution of $N$, which is locally free over $\Omega_Y^\gr$. We have the natural morphism: 
$$F\otimes_{\Omega_Y}f_*I\to f_\bullet(f^+F\otimes_{\Omega_X}I).$$

Take the composition $\Phi$ of this morphism with $$f_\bullet(f^+F\otimes_{\Omega_X}I)\to f_*(f^+F\otimes_{\Omega_X}I).$$

We claim that $\Phi$ is an isomorphism. It is enough to check the claim for graded modules. Clearly, $\Phi$ is an isomorphism when $F=\Omega_Y^\gr$. Thus the claim holds for any locally free $\Omega_Y^\gr$-module. The projection formula now follows from the definition of the (shifted) derived tensor product and of the extraordinary inverse image.
\qed

\subsection{Extraordinary inverse image for a closed embedding.}\label{inj_closed}
Let $i:Z\to X$ be a closed embedding of smooth schemes over $S$. In this case $\Omega_Z$ is coherent over $\Omega_X$, and we can give the second definition of the extraordinary inverse image:
$$i^b(M)=\cHom_{i^\bullet\Omega_X}(\Omega_Z,i^\bullet(M)).$$   

We define the derived functor $Ri^b:\D^\co(\qc{Y})\to\D^\co(\qc{X})$ using the equivalence~\ref{inj_resolution}.

Let $I$ be an injective $\Omega_X$-module. Since the functors $i_\bullet$ and $\cHom_{i^\bullet\Omega_X^\gr}(\Omega_Z^\gr,i^\bullet(-))$ are adjoint, and $i_\bullet$ is exact, 
the $\Omega_X$-module $i^b(I)$ is injective. It follows that for any pair of composable closed embeddings $f$ and $g$ we have $(f\circ g)^b\cong g^b\circ f^b$. For the same reason, we have:

\begin{thm*} The functors $i_*$ and $i^b$ are adjoint.
Let $N\in\D^\co(\qc{X})$ and let $M\in\D^\abs(\coh{X})$. Then
\begin{equation}\label{eq_cong_closed}
i_*\RHom_{\Omega_Z}(M,Ri^b(N)))\cong\RHom_{\Omega_X}(i_*M,N).
\end{equation} 
\end{thm*}
%\proof
%We prove (\ref{eq_cong_closed}). 
%The functors are adjoint by the argument from the proof of Theorem~\ref{adj_smooth}.
%Take an injective resolution $I$ for $N$. By Lemma $i^b(I)$ is injective. The result follows from the adjointness of $i_\bullet$ and $\cHom_{i^\bullet\Omega_X}(\Omega_Z,i^\bullet(-))$.
% By Theorem~\ref{thick_hom} we have to prove that
%$$i_*\cHom_{\Omega_Z}(P,\cHom_{i^\bullet\Omega_X}(\Omega_Z,i^\bullet(I)))\cong\cHom_{\Omega_X}(i_*P,I).$$  
%This is straightforward. 
% By Proposition~\ref{direct},
%$$\Hom_{\D^\co(\qc{Z})}(M,i^b(N)))\cong\Hom_{\D^\co(\qc{X})}(i_*M,N).\qed$$ 
%\qed

In~\ref{func_closed} we prove that if $S$ is defined over $\QQ$, then $Ri^b$ and $i^!$ are equivalent.

\begin{corollary}\label{cor_closed}
Let $i:Z\to X$ be a closed embedding, and $j:U=X\setminus Z\to X$ be the corresponding open embedding. Then for any $\Omega_U$-module $N$ we have: $i^bj_*N=0$.
\end{corollary}
\proof
By Corollary~\ref{inverse}.(3), for any $\Omega_Z$-module $M$ we have $j^+i_*M=0$.
The result follows from Theorems~\ref{adj_smooth} and~\ref{inj_closed}.
\qed

\section{Differential graded commutative algebra}\label{algebra}
\subsection{The injective hull.}\label{inj_hull}
Let $M$ be an $\Omega_X$-module. We say that injective $\Omega_X^\gr$-module $E$ is an \emph{injective hull of $M^\gr$}, if $M^\gr\subset E$, and every nonzero submodule of $E$ intersects $M^\gr$ non-trivially. 
The injective hull of a graded module exists and is unique up to an isomorphism~\cite{eisenbud}.
We say that injective $\Omega_X$-module $E$ is an \emph{injective hull of $M$} if the graded module $E^\gr$ is the injective hull of $M^\gr$.

\begin{prop*}
For any $\Omega_X$-module $M$ there exists an injective hull $E(M)$ of $M$.
\end{prop*}
\proof
Take an injective $\Omega_X$-module $I$ which contains $M$ as a submodule. Let $E^\gr$ be the injective hull of the graded module $M^\gr$ in $I^\gr$. The module $E^\gr$ generates an $\Omega_X$-submodule $E_1$ of $I$. Since $E^\gr$ is injective, $E_1^\gr\cong E^\gr\oplus E_2$ for some graded $\Omega_X^\gr$-submodule $E_2$ of $E_1^\gr$. The $\Omega_X^\gr$-module $E_2$ can be generated by (some set of) linear combinations of $de$ for $e\in E^\gr$, so $E_2$ is an $\Omega_X$-submodule of $E_1$. Put $E(M)=E_1/E_2$. 
\qed

\begin{rem*}
By definition, the graded module $E(M)^\gr$ is unique up to isomorphism. However, here is an example of an $\Omega_X$-module $M$ such that $E(M)^\gr$ has two different differentials compatible with differential on $M$. Suppose $X$ is affine.
Let $M$ be an $\Omega_X$-module such that $df\cdot M=0$ for some $f\in\OO_X(X)$. For example, $M=df\wedge\Omega_X$.
By the Proposition, there exists a module $E(M)$ with some differential $d$. Then $d+df$ is the second differential on $E(M)$,  and $dm=(d+df)m$ for any $m\in M(X)$.
\end{rem*}

\subsection{The stupid filtration.} Let $M$ be an $\Omega_X$-module. Define \emph{the stupid filtration $\st_\bullet$} on $M$ by the formula: $\st_i M=\Omega_X^iM$. Clearly, $\st_i M$ is an $\Omega_X$-submodule of $M$.
The associated graded module $\Gr^\st M$ is a DG-module over $\Gr^\st\Omega_X=\Omega_X^\gr$. In particular, the differential on $\Gr^\st M$ is $\Omega_X^\gr$-linear. We define the stupid filtration for $\Omega_X^\gr$-modules in the same way.

\begin{prop}\label{tor1}
Let $M$ be an $\Omega_X^\gr$-module. Suppose that $\cTor^1_{\Omega_X^\gr}(\OO_X,M)=0$. Then $\Gr_i^\st M\cong\Omega_X^i\otimes_{\Omega_X^\gr} M,$ and $\Gr^\st M\cong\Omega_X^\gr\otimes_{\OO_X}\Gr^\st_0 M$.
\end{prop}
\proof
Take two short exact sequences of perfect $\Omega_X^\gr$-modules:
$$0\to\st_i\Omega_X\to\Omega_X\to\Omega_X/\st_i\Omega_X\to 0, \text{ and}$$ 
$$0\to\st_{i+1}\Omega_X\to\st_i\Omega_X\to\Omega^i_X\to 0.$$ 
Take the tensor product over $\Omega_X^\gr$ of these sequences with $M$. We get that $\st_iM\cong \st_i\Omega_X\otimes_{\Omega_X^\gr} M$, and that $$\Gr^\st_i M\cong\Omega_X^i\otimes_{\Omega_X^\gr} M \cong \Omega_X^i\otimes_{\OO_X}\Gr^\st_0 M.\qed$$ 

%Since $\Omega^i$ is locally free, it follows that $\Gr^\st$ is exact on thick $\Omega_X^\gr$-modules.
\begin{corollary}\label{gr_st}
The functor $\Gr^\st$ is exact on the category of thick $\Omega_X^\gr$-modules.
\end{corollary}
\proof
The functors $\Gr^\st_0\cong\OO_X\otimes_{\Omega_X^\gr}(-)$ and $\Omega_X^i\otimes_{\OO_X}(-)$ are exact on $\thick{X}$.
\qed

\begin{thm}\label{thick_mod2}
 The following conditions on an $\Omega_X^\gr$-module $A$ are equivalent:
\begin{enumerate}
\item $A$ is thick;
\item $\cTor^1_{\Omega_X^\gr}(\OO_X,A)=0$;
\item $\Gr^\st A$ is (co)induced;
\item $\Gr^\st A$ is thick.
\end{enumerate}
\end{thm}
\proof
By Theorem~\ref{thick_mod}, $(1)\Rightarrow(2)$ and $(3)\Rightarrow(4)$. %$(3)$ and $(4)$.
By Proposition~\ref{tor1}, $(2)\Rightarrow(3)$. We now prove $(4)\Rightarrow(1)$.
Take a finite resolution for $\Gr^\st A$ by modules which are locally free over $\Omega_X^\gr$:
$$0\to P_n\to\dots\to P_0\to \Gr^\st A\to 0.$$ 
Let $P_0=Q_0\otimes_{\OO_X}\Omega_X^\gr$. Then there is a surjection $\phi:Q_0\cong\Gr^\st_0 P_0\to \Gr^\st_0 A$. For any point of $X$ there exists a neighborhood $U$ such that the restriction $\phi_{|U}$ lifts to a surjection of $\Omega^\gr_U$-modules $P^\gr_{0|U}\to A^\gr_{|U}$. This gives a surjection of $\Omega_U$-modules $\tilde\phi:\tilde P_0=\DG^+(Q_{0|U})\to A_{|U}$, and the projective dimension $$\mathrm{pd}_{\OO_U}\Gr^\st(\ker\tilde\phi)=\mathrm{pd}_{\OO_U}\Gr^\st A-1.$$
Similarly, we can make $U$ smaller and construct a finite resolution for $A_{|U}$ by thick $\Omega_U$-modules
$$0\to \tilde P_n\to\dots\to \tilde P_0\to A_{|U}\to 0.$$ 
%such that functor $\Gr^\st$ takes this resolution to $(*)$.
It follows that $A_{|U}$ is thick. Thus for any $i>0$ the sheaf $\cTor^i_{\Omega_X^\gr}(\OO_X,A)$ is trivial.
\qed

\begin{ex} Here is an example of a thick $\Omega_X^\gr$-module $M$ which is not induced. Let $X=\Spec k[t]$, and let $\delta=k[t]/tk[t]$ be a skyscraper module. Put $M=k[t]\oplus k[t]\frac{dt}{t}\oplus\delta dt$, where $k[t]\frac{dt}{t}$ is the submodule of $k(t)dt$ generated over $k[t]$ by $\frac{dt}{t}$.
The action of $t$ is natural, and to define the action of $dt$ one should write $k[t]\frac{dt}{t}$ as the sum $k[t]dt\oplus\delta$. We see that $Gr^\st M\cong\Omega_X^\gr\oplus\delta\otimes_{\OO_X}\Omega_X^\gr$. Thus $M$ is thick and is not induced.
\end{ex}

\section{Kashivara theorem for $\Omega_X$-modules}\label{kt}
\subsection{}\label{MV}
For an open subset $U\subset X$ we denote by $j_U:U\to X$ the corresponding open embedding. Recall that $j_U^!\cong j_U^\bullet(-)\otimes_{j_U^\bullet\OO_X}\OO_U$. This justifies the following notation:  $M_{|U}=j_U^!M$, where $M$ is an $\Omega_X$-module.

\begin{lemma*}
Let $M$ be an $\Omega_X$-module, and let $U$ and $V$ be open subsets of $X$. Then we have the Mayer-Vietoris triangle in $\D^\co(\qc{X})$: 
$$j_{U\cup V*}j_{U\cup V}^!M\to j_{U*}j_{U}^!M\oplus j_{V*}j_{V}^!M\to j_{U\cap V*}j_{U\cap V}^!M\to.$$
\end{lemma*}
\proof
Let $I$ be an injective representative for $M$. 
Then the sequence of $\Omega_X$-modules:
$$0\to j_{U\cup V*}I_{|U\cup V}\to j_{U*}I_{|U}\oplus j_{V*}I_{|V}\to j_{U\cap V*}I_{|U\cap V}\to 0$$
is exact since $I$ is flabby.
\qed

It follows that coacyclicity is a local property.
\begin{corollary*}
Let $M$ be an $\Omega_X$-module. Suppose there is a finite covering $\mathcal{U}$ of $X$ such that for any $U\in \mathcal{U}$ the restriction $M_{|U}$ is coacyclic. Then $M$ is coacyclic.
\end{corollary*}

\subsection{} Let $i:Z\to X$ be a closed embedding, % of smooth schemes over $S$,
 and $j:U=X\setminus Z\to X$ be the corresponding open embedding. Define the functor of ``sections with support in $Z$'' 
$$\Gamma^X_Z:\qc{X}\to\qc{X}$$ as follows. Let $J_Z\subset\OO_X$ be the ideal sheaf of $Z$. Then the DG-ideal sheaf of $Z$ is the DG-ideal $J\subset\Omega_X$ generated by $J_Z$ and $dJ_Z$.
 For an $\Omega_X$-module $M$ put
$$\Gamma^X_Z(M)=\varinjlim_{m}\cHom_{\Omega_X}(\Omega_X/J^m, M).$$
In fact, $\Gamma^X_Z(M)\cong \varinjlim_{m}\cHom_{\OO_X}(\OO_X/J_Z^m, M).$
Using the argument similar to~\cite[III.3.2]{Ha} one proves that for an injective $\Omega_X$-module $I$ the module $\Gamma^X_Z(I)$ is also injective.
Clearly, $\Gamma^X_Z$ is exact on injective $\Omega_X$-modules. Thus we can define the derived functor $R\Gamma^X_Z:\D^\co(\qc{X})\to\D^\co(\qc{X})$ using injective resolutions. When $X$ is fixed we often write $\Gamma_Z$ for $\Gamma^X_Z$. 

\begin{prop*}\label{first_triangle}
For any $\Omega_X$-module $M$ there is a functorial distinguished triangle in $\D^\co(\qc{X})$:
$$R\Gamma_Z(M)\to M\to j_*j^!M\to$$ 
\end{prop*}
\proof
We have to prove that for an injective $\Omega_X$-module $I$ there is an exact triple: $$0\to\Gamma_Z(I)\to I\to j_*j^!I\to 0.$$ The left morphism is the natural inclusion, and the right morphism comes from~\ref{adj_smooth}. Since $I$ is injective, the triple is left and right exact. Now we prove that the triple is exact in the middle. We may assume that $X$ is affine. Let $m\in I(X)$. If the restriction of $m$ on $U$ is zero, then for any $f\in J_Z$ there exist $k(f)$ such that $f^{k(f)}m=0$. Since $J_Z$ is finitely generated there exists $k$ such that $J^km=0$, i.e. $m\in\Gamma_Z(I)$.
\qed

\begin{thm}\label{pre_Kashivara} Assume that $S$ is defined over $\QQ$, i.e. all natural numbers are invertible in $S$. If $Z$ is smooth, then for any $\Omega_X$-module $M$ we have $i_*Ri^b(M)\cong R\Gamma_Z(M)$ in $\D^\co(\qc{X})$. 
\end{thm}
\proof
Take an injective resolution $I$ of $M$. Then $i_*i^b(I)\cong\cHom_{\Omega_X}(i_*\Omega_Z,I)$ is a submodule of $\Gamma_Z(I)$. This gives a morphism of functors. It is enough to prove that this morphism is equivalence when $X$ is affine, and $J_Z$ is generated by a regular sequence $f_1,\dots, f_r$. 

\begin{sublemma}
Let $a:Z_1\to X$ and $Z_2\subset Z_1$ be closed embeddings. Then
$$R\Gamma^{Z_1}_{Z_2}\circ Ra^b\cong Ra^b\circ R\Gamma^X_{Z_2}.$$
\end{sublemma}
\proof
Apply both sides to $I$.
The module $E=\Gamma^X_{Z_2}(I)$ is injective, thus $a^b(E)$ computes the right hand side. By Proposition~\ref{first_triangle}, we have the short exact sequence of $\Omega_X$-modules: $$0\to E\to I\to j_*j^!I\to 0,$$ where $j:X\setminus Z_2\to X$ is the natural open embedding. Apply the functor $\Gamma^{Z_1}_{Z_2}\circ a^b$ to this exact sequence. By Corollary~\ref{cor_closed}, $a^bj_*j^!I\cong 0$, and it is easy to check that
$\Gamma^{Z_1}_{Z_2}\circ a^b(E)\cong a^b(E)$.
The lemma follows. 
\qed

By the Lemma, we may assume that the ideal of $Z$ in $\OO_X$ is generated by a single regular function $f$.

Let $N=\Gamma_Z(I)$, and let $N_0\cong i_*i^b(I)$ be the submodule of $N$ annihilated by $f$ and $df$. Put $$N_k=\{m\in N| f^km\in N_0\}.$$ By Lemma~\ref{lem1} below, $N_k$ is an $\Omega_X$-module.
Clearly, $N\cong\varinjlim_k N_k$.
We claim that for $k>0$ the module $N_k/N_{k-1}$ is contractible. Then by Lemma~\ref{colemma}(3), $N$ is co-isomorphic to $N_0$. 

We prove the claim. 

\begin{sublemma}\label{lem1}
If $m\in N_k$, then $dm\in N_k$.
\end{sublemma}
\proof
We prove that $f^kdm\in N_0$, i.e. $f(f^kdm)=0$, and $df(f^kdm)=0$. This is clear from the relations:
$$0=d(f^{k+1}m)=(k+1)f^kdf\cdot m+f^{k+1}dm,\text{ and}$$ 
$$0=d(f^kdf\cdot m)=-f^kdf\cdot dm.$$
\qed

%\begin{sublemma}\label{lem2}
%If $m\in N_k$, and $k$ is invertible on $S$, then $df\cdot m\in N_{k-1}$.
%\end{sublemma}
%\proof Indeed, $$d(f^km)=kf^{k-1}df\cdot m+f^kdm$$ is an element of $N_0$, and $f^kdm\in N_0$. It follows that $f^{k-1}(df\cdot m)\in N_0$.\qed

\begin{sublemma}\label{lem4}
If $df\cdot m=0$ for some $m\in I$, then there exists $\beta\in I$ such that $m=df\cdot\beta$.
\end{sublemma}
\proof
Recall that $f$ is regular, and $I$ is injective. Define the morphism $\phi: df\wedge\Omega_X\to I$ by $\phi(df\wedge\omega)=\omega\cdot m$. This definition is correct since $f$ is regular. We can lift $\phi$ to a morphism $\tilde\phi:\Omega_X\to I$. Put $\beta=\tilde\phi(1)$.
\qed

\begin{sublemma}\label{lem3}
Suppose that $m\in N_k$. Then there exists $\alpha\in N_k$ such that $f^km=f^{k-1}df\cdot\alpha$, and $f^k\alpha=0$.
\end{sublemma}
\proof
Put $m'=f^km\in N_0$. Then $df\cdot m'=0$, and, by Lemma~\ref{lem4}, there exists $\beta\in I$ such that $m'=df\cdot\beta$. By assumption, $$fdf\cdot\beta=fm'=0.$$ As before, there exists $\gamma\in I$ such that $f\beta=fdf\cdot\gamma$. Put $\alpha'=\beta-df\cdot\gamma$. Then  $f\alpha'=0$, and $m'=df\cdot\alpha'$. There exists $\alpha\in I$ such that $f^{k-1}\alpha=\alpha'$. Thus $f^k\alpha=0$, and $\alpha\in N_k$. 
\qed

We define an endomorphism $s$ on $N_k/N_{k-1}$ as follows. Take $m\in N_k$. Then, by Lemma~\ref{lem3}, there exists $\alpha\in N_k$ such that $f^km=f^{k-1}df\cdot\alpha$. Let $s(m)$ be the image of $\alpha$ in $N_k/N_{k-1}$. Clearly, $s(m)$ is independent on the choice of $\alpha$, and 
%Since $f^k\alpha=0$, we have $s^2=0$. By Lemma~\ref{lem2}, $df$ annihilates $N_k/N_{k-1}$, thus 
$s$ is a morphism of $\Omega_X^\gr$-modules of degree $-1$. Let us prove that $ds+sd=-k$. From the relation $f^km=f^{k-1}df\cdot \alpha$ we get:
$$f^kdm=f^{k-1}df(-d\alpha-km).$$ It follows that $s(dm)=-ds(m)-km$, and $-s/k$ is a contracting homotopy for $N_k/N_{k-1}$.
\qed

\begin{corollary}\label{cor_Kashivara}
Assume that $S$ is defined over $\QQ$. Let $I$ be an injective $\Omega_Z$-module. Then an injective hull $E$ of $i_*I$ is co-isomorphic to $i_*I$.
\end{corollary}
\proof
Note that $E=\Gamma_Z(E)$, because it is true for the graded module $E^\gr$. By Theorem~\ref{pre_Kashivara}, $E\cong i_*i^b(E)$ in $\D^\co(\qc{X})$. We prove that $i^b(E)$ is isomorphic to $I$ over $\Omega_Z$. Since $I$ is injective, it is a direct summand of $i^b(E)$ as a graded module, i.e. $i^b(E)^\gr\cong I^\gr\oplus I'$ for some $\Omega_Z^\gr$-module $I'$. Thus $i_*I'$ is a submodule of $E^\gr$ such that $i_\bullet I'\cap i_\bullet I=0$. Since $E^\gr$ is an injective hull of $i_\bullet I^\gr$, we conclude that $i_\bullet I'=0$. The functor $i_\bullet$ is full and faithful, thus $I'=0$.
\qed

\begin{thm}{\rm(Kashivara theorem.)}\label{Kashivara}
Assume that $S$ is defined over $\QQ$.
Let $\D_Z^\co(\qc{X})$ be the full subcategory of $\D^\co(\qc{X})$ of $\Omega_X$-modules with support on $Z$, and $\D_Z^\abs(\coh{X})$ be the corresponding category of coherent $\Omega_X$-modules. Then the functor $i_*$ gives an equivalence of $\D^\co(\qc{Z})$  with $\D_Z^\co(\qc{X})$, and of $\D^\abs(\coh{Z})$ with $\D_Z^\abs(\coh{X})$. In particular, $i_*$ is full and faithful.
\end{thm}
\proof%[Proof of Kashivara theorem.]
Let $M$ be an object of $\D_Z^\co(\qc{X})$.
By Theorem~\ref{pre_Kashivara} and Proposition~\ref{first_triangle}, we have the distinguished triangle in $\D^\co(\qc{X})$:
$$i_*i^b(M)\to M\to j_*j^+M\to.$$ By assumption, $j^+M=0$, thus $i_*i^b\cong Id_{\D_Z^\co(\qc{X})},$ and
by Corollary~\ref{cor_Kashivara}, $i^bi_*\cong Id_{\D^\co(\qc{Z})}.$

For a coherent $\Omega_Z$-module $M$ the module $i_*M$ is also coherent. Conversely, any coherent $\Omega_X$-module $M$ with support in $Z$ has a finite filtration such that the associated graded module is the direct image $i_*M'$ of a coherent $\Omega_Z$-module $M'$. This is enough, since $\D^\abs(\coh{X})$ is a full subcategory of $\D^\co(\qc{X})$~\cite[3.11]{P2}.
\qed

%\begin{sublemma}
%For any $\Omega_X$-module $M$ we have $i^bj_*M=0$.
%\end{sublemma}

\begin{thm}{\rm(Base change.)} Assume that $S$ is defined over $\QQ$.
Let $f:X\to Y$ and $g:Z\to Y$ be two morphisms of smooth schemes over $S$. Consider the Cartesian square
$$\begin{CD}
W@>\tilde{g}>>X\\
@V\tilde{f} VV @V f VV\\
Z@> g>> Y
\end{CD}$$
Assume that $W=X\times_Y Z$ is smooth over $S$. Then there exists an isomorphism 
$$g^!\circ f_*\cong\tilde f_*\circ\tilde g^!:\D^\co(\qc{X})\to\D^\co(\qc{Z}).$$
\end{thm}
The proof is the same as in~\cite[1.7.3]{HTT}.
\qed

\section{Trace map for a closed embedding}\label{tce}
\subsection{\v{C}ech resolution.}\label{cech}
Let $\UU=\{U_i\}_{i\in I}$ be a finite open covering of $X$, where $I=\{0,\dots, N\}$ is a finite ordered set. 
Let $M$ be an $\Omega_X$-module. We define the \v{C}ech resolution $C^\bullet(\UU,M)$ of $M$ as follows: $$C^n(\UU,M)=\oplus_{J\subset I}j_{J*}j_J^+M,$$ where $J$ is a subset of $I$ of order $|J|=n+1$, and $j_J$ is the  inclusion $\cap_{i\in J}U_i\to X$. The closed morphism $\d_n:C^n(\UU,M)\to C^{n+1}(\UU,M)$ is defined 
%as follows. Let $m=(m_J)_{J\subset I}\in C^n(\UU,M)$, where $m_I\in j_{J*}j_J^+M$. Then $\d_n(m)_{J'}=\sum_{J'=J\cup\{i\}}(-1)^im_J$. Clearly $\d_n$ is closed.
in the usual way. 
Put $C(\UU,M)=\Tot^\oplus C^\bullet(\UU,M)$. The natural closed morphisms $M\to j_{i*}M_{|U_i}$ induce the natural inclusion morphism $M\to C^0(\UU,M)\to C(\UU,M)$. 

\begin{lemma*}
The closed morphism $M\to C(\UU,M)$ is a co-isomorphism. Moreover, its cone is absolutely acyclic.
\end{lemma*}
\proof
The sequence of $\Omega_X^\gr$-modules
$$0\to M^\gr\to C^0(\UU,M)^\gr\to\dots \to C^{|I|}(\UU,M)^\gr\to 0$$ 
is exact, because it is the \v{C}ech resolution of the graded quasi-coherent $\OO_X$-module $M^\gr$. 
The $\Omega_X$-module $\cone(M\to C(\UU,M))$ is the total object of the sequence
$$0\to M\to C^0(\UU,M)\to\dots \to C^{|I|}(\UU,M)\to 0.$$ It is absolutely acyclic, since we can cut this sequence into short exact sequences. The lemma follows.
\qed

\begin{thm*}
Suppose $\mathcal{U}$ is a finite affine covering of $X$. Then $f_*(M)$ is co-isomorphic to $f_\bullet(C(\UU,M))$.
\end{thm*}
\proof
Let $N$ be an injective resolution of $M$. Then $j_J^+(N)$ is flabby for any $J\subset I$, thus $C(\UU,N)$ is flabby too, and $f_*(M)$ is co-isomorphic to $f_\bullet(N)$. By the Lemma, $f_\bullet(N)$ is co-isomorphic to $f_\bullet C(\UU,N)$.
There is the natural morphism $\phi:C(\UU,M)\to C(\UU,N)$ with coacyclic cone.
Since the morphism $g=f\circ j_J$ is affine, the functor $g_\bullet$ is exact, and $g_\bullet j_J^+(N)$ is co-isomorphic to $g_\bullet j_J^+(M)$. 
%By definition, $C(M)=\Tot C^\bullet(M)$. This gives a stupid filtration on $C(M)$ and a compatible stupid %filtration on $C(N)$. We proved that associated graded object with respect to this filtrations  are co-isomorphic. 
It follows that $\phi$ induces a morphism $f_\bullet C(\UU,M)\to f_\bullet C(\UU,N)$ with coacyclic cone.  
\qed

\subsection{Gluing morphisms in $\D^\co(\qc{X})$.}
In this section we prove a counterpart of \cite[3.2.2]{BBD} for $\Omega_X$-modules. 

\begin{thm*}\label{glue}
Let $M$ and $N$ be $\Omega_X$-modules. Suppose that $M$ is coherent.
Let $\mathcal{U}$ be a finite affine covering of $X$. Assume that $\Hom_{\D^\co(\qc{V})}(M_{|V},N_{|V}[j])=0$ for all affine $V\subset U\in\mathcal{U} $ and $j<0$. Then $$\cHom_{\D^\co(\qc{X})}(M,N):V\mapsto \Hom_{\D^\co(\qc{V})}(M_{|V},N_{|V})$$ is a sheaf.
\end{thm*}
\proof
It is enough to prove that sections over $X$ can be glued of sections over the covering $\mathcal{U}$. We may assume that $N$ is injective. %For an open subscheme $U$ the module $\cHom_{\Omega_U}(M_{|U},N_{|U})$ can be viewed as a complex of $\OO_S$-modules. 
Define the presheaf $$\mathcal{H}^q:V\mapsto\Hom_{\D^\co(\qc{V})}(M_{|V},N_{|V}[q]).$$ 
We have to prove that $\mathcal{H}^0$ is a sheaf. By proposition~\ref{direct},
$$\mathcal{H}^q(V)\cong H^q(S,f_{|V*}\cHom_{\Omega_V}(M_{|V},N_{|V})),$$ where $f:X\to S$ is the morphism to the base scheme, and $f_{|V}$ is the restriction to $V$.
By Theorem~\ref{cech}, we have the \v Cech spectral sequence: $$H^p(\mathcal{U},\mathcal{H}^q))=>H^{p+q}(S,f_*\cHom_{\Omega_X}(M,N))\cong \Hom_{\D^\co(\qc{X})}(M,N[p+q]).$$
%\Hom_{\D^\co(\qc{X})}(M,N[p+q]).$$
%By Proposition~\ref{direct}, $\Hom_{\D^\co(\qc{U})}(M_{|U},N_{|U})$ is canonically isomorphic to $H^0(S,\Hom_{\Omega_U}(M,N))$; thus there is an isomorphism of presheaves $\cHom_{\D^\co(\qc{X})}(M,N[q])\cong\mathcal{H}^q$. 
%Since $U$ is affine, $\mathcal{H}^q_{|U}$ is a quasi-coherent sheaf.
By assumption, $H^p(\mathcal{U},\mathcal{H}^q))=0$ for $q<0$. It follows that $H^0(\mathcal{U},\mathcal{H}^0)\cong \Hom_{\D^\co(\qc{X})}(M,N)$. 
\qed

\subsection{Koszul resolution.}\label{koszul}
Let $i:Z\to X$ be a closed embedding of affine smooth schemes over $S$. Suppose that $Z$ is a complete intersection in $X$. In this section we construct an explicit perfect resolution for $\Omega_Z$ over $\Omega_X$. 

\begin{thm*}
There exists a perfect $\Omega_X$-module $K$ and a non-canonical co-isomorphism $K[n]\to i_*\Omega_Z$.
\end{thm*}
\begin{sublemma}
Suppose that the ideal of $Z$ in $\OO_X$ is generated by $f_1,\dots,f_n\in\OO_X$, where $n=\dim X-\dim Z$.
Let $P_Z$ be the submodule of $\Omega_X$ generated by $df_1\wedge\dots\wedge df_n$, and let $I_Z$ be the submodule $f_1P_Z+\dots+f_nP_Z$. Then $P_Z$ is perfect, and $i_*\Omega_Z\cong P_Z/I_Z[n]$.
\end{sublemma}\label{lem31}
\proof
Locally there exist functions $z_1,\dots,z_s\in\OO_X$ such that $dz_1,\dots,dz_s,df_1\dots df_n$ generate $\Omega^1_X$ over $\OO_X$. Thus $P_Z$ is a direct summand in $\Omega_X$ and perfect. We have the natural surjective morphism $P_Z[n]\to i_*\Omega_Z$ given by $df_1\wedge\dots\wedge df_n\mapsto 1$. Clearly, its kernel is $I_Z$.
\qed

\begin{sublemma}\label{lem32}
Let $L_\bullet$ be the Koszul resolution of $\OO_Z$ over $\OO_X$. Take the complex of $\Omega^\gr_X$-modules: $K_\bullet=L_\bullet\otimes_{\OO_X}P_Z$. Then the induced morphism $\delta_{m}:K_m\to K_{m-1}$ is a closed morphism of $\Omega_X$-modules.  
\end{sublemma}
\proof
By definition of the Koszul complex, $L_\bullet\cong\wedge^\bullet L_1$ is a DG-algebra over $\OO_X$ generated by $L_1$, where $L_1\cong\oplus_{i=1}^n v_i\OO_X$ for some $v_i\in L_1$ and the differential is given by $v_i\mapsto f_i$.
Thus  $K_0\cong P_Z$, and $K_1\cong\oplus_{i=1}^n v_iP_Z$.
The differential $\delta_{m}$ on $K_\bullet$ comes from the formula: $$\delta_1(v_i\otimes df_1\wedge\dots\wedge df_n)=f_idf_1\wedge\dots\wedge df_n.$$ Now it is clear that $\delta_{m}$ is a closed morphism of $\Omega_X$-modules.
\qed

\proof[Proof of Theorem~\ref{koszul}]
Clearly, the sequence of $\Omega_X$-modules
$$0\to K_n\to\dots\to K_1\to K_0\cong P_Z\to i_*\Omega_Z[-n]\to 0$$ is exact.
Put $K=\Tot^\oplus(K_\bullet)$.
\qed

\begin{thm}\label{ftrace} Suppose that $S$ is a $\QQ$-scheme.
Then there is a natural co-isomorphism $\tr:i^b(\Omega_X)\cong \Omega_Z[2\dim X-2\dim Z]$.
\end{thm}
\proof
Apply $i_*$ to both sides.
By Theorem~\ref{Kashivara}, it is enough to prove that $$\RHom_{\Omega_X}(i_*\Omega_Z,\Omega_X)\cong i_*\Omega_Z[2\dim X-2\dim Z].$$
We construct the isomorphism locally and then use Theorem~\ref{glue}.

Assume that $X$ is affine and the ideal of $Z$ in $\OO_X$ is generated by $f_1,\dots,f_n\in\OO_X$, where $n=\dim X-\dim Z$. Let $K[n]$ be the Koszul resolution for $i_*\Omega_Z$. By Theorem~\ref{thick_hom}, $$\RHom_{\Omega_X}(K[n],\Omega_X)\cong\cHom_{\Omega_X}(K[n],\Omega_X).$$ 
We use the notation from the proof of Lemma~\ref{lem32}. By definition, $K_1\cong\oplus_{i=1}^n v_iP_Z$. Hence $K_n\cong P_Z$ is generated by $v_1\wedge\dots\wedge v_n$. Moreover, $\cHom_{\Omega_X}(P_Z,\Omega_X)$ is isomorphic to $P_Z[n]$ non-canonically. Thus the complex of $\Omega_X$-modules $\Check K^\bullet=\cHom_{\Omega_X}(K_\bullet[n],\Omega_X)$ is isomorphic to $K_\bullet$. 
The module $\Check K^n[n]=\cHom_{\Omega_X}(K_n,\Omega_X)\cong P_Z$ has the generator $u$ given by $$v_1\wedge\dots\wedge v_n\mapsto df_1\wedge\dots\wedge df_n.$$
As in Lemma~\ref{lem31}, we define the morphism $\Check K^n[n]\to \Omega_Z$ such that $u$ goes to $1\in\Omega_Z$. This gives a co-isomorphism $$\psi_X:\Check K=\RHom_{\Omega_X}(i_*\Omega_Z,\Omega_X)\to i_*\Omega_Z[-2n].$$ It is clear that $\psi_X$ does not depend on the choice of coordinates $f_1,\dots,f_n$.

We proved the theorem locally and now turn to the global situation. We have to find a co-isomorphism of $\Check K=\RHom_{\Omega_X}(i_*\Omega_Z,\Omega_X)$ and $i_*\Omega_Z[-2n]$. By the local result,
\begin{equation*}
\begin{split}
&\Hom_{\D^\co(\qc{V})}(\Check K_{|V},i_*\Omega_{Z|V}[-2n+j])\cong
\Hom_{\D^\co(\qc{V})}(i_*\Omega_{Z|V},i_*\Omega_{Z|V}[j])\cong\\
&\Hom_{\D^\co(\qc{Z\cap V})}(\Omega_{Z\cap V},\Omega_{Z\cap V}[j])=0
\end{split}
\end{equation*}
for small affine open $V\subset X$ and any $j<0$.
By Theorem~\ref{glue}, the presheaf $\cHom_{\D^\co(\qc{X})}(\Check K,i_*\Omega_Z[-2n])$ is a sheaf. Thus we can glue local morphisms $\psi_U$ to a global morphism $\psi\in\Hom_{\D^\co(\qc{X})}(\Check K,i_*\Omega_Z[-2n])$. This morphism is a co-isomorphism by Corollary~\ref{MV}.
\qed

\begin{thm}\label{func_closed} Suppose that $S$ is a $\QQ$-scheme.
Let $i:Z\to X$ be a closed embedding of smooth schemes over $S$. 
Then the functors $i^b$ and $i^!$ are naturally equivalent.
\end{thm}
\proof
In Theorem~\ref{ftrace} we construct the trace isomorphism for closed embedding: $$\tr:\RHom_{f^\bullet\Omega_X}(\Omega_Z,f^\bullet\Omega_X)\cong \Omega_Z[2\dim Z-2\dim X].$$
By Remark~\ref{perf_resolution}, any $\Omega_X$-module has a resolution $F$ which is
locally free over $\Omega^\gr_X$. We have 
\begin{equation*}
\begin{split}
i^b(F)=&\RHom_{f^\bullet\Omega_X}(\Omega_Z,f^\bullet(F))\cong  \RHom_{f^\bullet\Omega_X}(\Omega_Z,f^\bullet\Omega_X)\otimes_{f^\bullet\Omega_X}f^\bullet(F)\cong\\ &\Omega_Z\otimes_{f^\bullet\Omega_X}f^\bullet(F)[2\dim Z-2\dim X]=i^!(F).\qed
\end{split}
\end{equation*}

\section{Coherent $\Omega$-modules}\label{sec_coh}

\subsection{} Recall that $\D^\abs(\coh{X})$ is a full subcategory of $\D^\co(\qc{X})$.
Moreover, coherent $\Omega_X$-modules are compact objects of $\D^\co(\qc{X})$~\cite[3.11]{P2}. Let $\D^\scoh(X)$ be the thick envelope of $\D^\abs(\coh{X})$ in $\D^\co(\qc{X})$.
We prove that an object $M$ of $\D^\co(\qc{X})$ belongs to $\D^\scoh(X)$ if and only if $M$ is locally coherent.

\begin{prop*}\cite[Remark 1.10]{P3}
Let $\mathcal{U}$ be a finite open covering of $X$, and let $M$ be an object of $\D^\abs(\coh{X})$. Then $M$ belongs to $\D^\scoh(X)$ if and only if for any $U\in\mathcal{U}$ the module $M_{|U}$ belongs to $\D^\scoh(U)$.
\end{prop*}
\proof We give a proof under assumption that $S$ is defined over $\QQ$.
Take $U\in\mathcal{U}$. By Kashivara theorem, $\D^\co(\qc{U})$ is equivalent to the Verdier quotient of $\D^\co(\qc{X})$ by $\D^\co(\qc{Z})$. Now, by Theorem~\cite[2.1.5]{Ne}, there is an object $N$ of $\D^\scoh(X)$ and a morphism $N\to M$ such that $N_{|U}\to M_{|U}$ is a co-isomorphism. The cone $M_1=\cone (N\to M)$ is coacyclic on $U$. Take  $V\in\mathcal{U}$, then $M_{1|V}$ belongs to $\D^\scoh(V)$. We can make $V$ smaller in such a way that $Z=V\setminus(U\cap V)$ is smooth. By Kashivara theorem, $M_{1|V\cup U}$ is co-isomorphic to $i_*M_Z$, where $i:Z\to V\cup U$ is the natural inclusion, and $M_Z$ is coherent. We proved that $M_{|U\cup V}$ belongs to $\D^\scoh(U\cup V)$. The proposition follows by induction.
\qed

\begin{thm}\label{coherence}
Let $f:X\to Y$ be a morphism of smooth schemes over $S$. If $f$ is proper, and $M$ is a coherent $\Omega_X$-module, then $f_*(M)$ is isomorphic to a coherent $\Omega_Y$-module in $\D^\co(\qc{Y})$.
If $f$ is smooth, and $N$ is a coherent $\Omega_Y$-module, then $f^!(N)$ is isomorphic to a coherent $\Omega_X$-module in $\D^\co(\qc{X})$.
\end{thm}
\proof
The category of coherent $\Omega_X$-modules is generated by coherent modules annihilated by $\Omega_X^1$. Suppose $M$ is such a module. We have to prove that $f_*(M)$ is coherent. The direct image can be computed using \v{C}ech resolution $C^\bullet(M)$. Clearly, $f_\bullet C^\bullet(M)$ is annihilated by $\Omega_Y^1$, and by Serre's theorem is quasi-isomorphic to a complex of coherent $\OO_Y$-modules. This quasi-isomorphism is naturally a co-isomorphism of 
$\Omega_Y$-modules.
The second part of the theorem follows from Lemma~\ref{inverse}. 
\qed

\subsection{Duality.}
Let $n=\dim_SX=\dim X-\dim S$ be a relative dimension of $X$ over $S$. We define the duality functor 
$$\D_X:\D^\abs(\coh{X})\to\D^\abs(\coh{X})$$ by $\D_X=\RHom_{\Omega_X}(-,\Omega_X[2n])$. Take a perfect $\Omega_X$-module $P$. Since $\Omega_X$ is coinduced the coherent module $\cHom_{\Omega_X}(P,\Omega_X[2n])$ represents $\RHom_{\Omega_X}(P,\Omega_X[2n])$. It follows that $\D_X$ is well-defined.

\begin{prop*}
$\D_X^2\cong \mathrm{Id}_{\D^\abs(\coh{X})}$.
\end{prop*}
\proof
It is enough to prove that $\D_X^2$ is equivalent to the identity functor on the category of prefect $\Omega_X$-modules. For a perfect $\Omega_X$-module $P$ we have a natural transformation $$\Phi:P\to \cHom_{\Omega_X}(\cHom_{\Omega_X}(P,\Omega_X),\Omega_X)\cong  \cHom_{\Omega_X}(\cHom_{\Omega_X}(P,\Omega_X[2n]),\Omega_X[2n]).$$
% given by $m\mapsto(\phi\mapsto\phi(m))$, where $m\in P$ and $\phi\in\cHom_{\Omega_X}(P,\Omega_X[2n])$. 
In fact, $\Phi$ is an isomorphism of $\Omega_X$-modules. We prove it locally for $P=\OO_X$. In this case $\cHom_{\Omega_X}(\OO_X,\Omega_X)\cong \Omega^n_X[-n]$, and $\cHom_{\Omega_X}(\Omega^n_X[-n],\Omega_X)\cong\OO_X$.
\qed

\begin{prop}\label{dual_hom}
Let $M$ and $N$ be coherent $\Omega_X$-modules. Then $D(M)\otimes^D_{\Omega_X}N\cong\RHom_{\Omega_X}(M,N)$.
% in $\D^\co(\qc{X})$.
\end{prop}
\proof
Choose a perfect resolution $P$ for $M$ and a coinduced resolution $A$ for $N$, i.e. 
$A^\gr\cong\cHom_{\OO_X}(\Omega_X^\gr,L)$ for some graded $\OO_X$-module $L$.
Let $\Phi$ be the natural transformation $$\cHom_{\Omega_X}(P,\Omega_X)\otimes_{\Omega_X}A\to\cHom_{\Omega_X}(P,A)$$
given by the formula: $\phi\otimes a\mapsto (p\mapsto\phi(p)a)$, where $\phi\in \cHom_{\Omega_X}(P,\Omega_X)$, $p\in P$, and $a\in A$.
We claim that $\Phi$ is an isomorphism of $\Omega_X$-modules. Indeed, the question is local and we may assume that $P=\OO_X$. In this case, $\Phi$ is an isomorphism since $A$ is coinduced.
The proposition follows.
\qed

\begin{prop}\label{dual_direct_closed_imb} Suppose $S$ is a $\QQ$-scheme.
Let $i:Z\to X$ be a closed embedding of smooth schemes over $S$. Then $\D_X\circ i_*\cong i_*\circ\D_Z$.
\end{prop}
\proof
By Theorems~\ref{inj_closed}~and~\ref{func_closed}, for any coherent $\Omega_Z$-module $M$ we have:
\begin{equation*}\begin{split}
\D_X(i_*M)=\RHom_{\Omega_X}(i_*M,\Omega_X[2\dim_S X])\cong i_*\RHom_{\Omega_Z}(M,i^!\Omega_X[2\dim_S X]))\cong \\ i_*\RHom_{\Omega_Z}(M,\Omega_Z[2\dim_S Z]))=i_*D_Z(M).
\end{split}\end{equation*}
\qed

\subsection{Trace for a proper projection.} Suppose that the base scheme $S$ is defined over $\QQ$. Let $p:\PP^n=\PP^n_S\to S$ be the natural projection from the projective space to $S$. We define the trace morphism $Tr_p$ as the composition $$f_*\Omega_{\PP^n}[2n]\to H^n(\PP^n,\Omega_{\PP^n}^n)\to\OO_S$$ of the natural projection to $H^n(\PP^n,\Omega_{\PP^n}^n)$ with the trace morphism defined in~\cite[Remark 10.2]{Ha}.

Let $g:Y\to S$ be a smooth $S$-scheme. Denote by $f:X=\PP^n\times_S Y\to Y$ and $q:X\to \PP^n$ natural projections to $Y$ and $\PP^n$. By the base change formula we have the trace morphism:
\begin{equation*}\begin{split}
Tr_f:f_*\Omega_X[2\dim Y+2n]\cong f_*(q^!\Omega_{\PP^n}[2n])\cong g^!p_*\Omega_{\PP^n}[2n]\xrightarrow{Tr_p} g^!\OO_S\cong \Omega_Y.
\end{split}\end{equation*}

\begin{lemma}
Suppose $f:X\to Y$ is smooth and proper. For any $M,N\in\D^\abs(\coh{X})$ there is a natural morphism $$f_*\RHom_{\Omega_X}(M,N)\to \RHom_{\Omega_X}(f_*M,f_*N).$$
\end{lemma}
\proof
By Theorem~\ref{coherence}, $f_*M$ is coherent, and the right hand side is well-defined.
By Theorem~\ref{adj_smooth}, we have an adjunction morphism $f^+f_*M\to M$, and the desired morphism is the composition: 
$$f_*\RHom_{\Omega_X}(M,N)\to f_*\RHom_{\Omega_X}(f^+f_*M,N)\to \RHom_{\Omega_X}(f_*M,f_*N).$$
\qed

\begin{thm}\label{dual_direct}
Let $f:X\to Y$ be a proper morphism of smooth schemes over $S$. Then $\D_Y\circ f_*\cong f_*\circ\D_X$.
\end{thm}
\proof
The morphism $f$ is a composition of a closed embedding and a proper projection. 
By Proposition~\ref{dual_direct_closed_imb}, we have to prove the theorem for the projection $f:X=\PP^n\times_S Y\to Y$. For any coherent $\Omega_X$-module $M$ we have the natural morphism:
\begin{equation*}\begin{split}
f_*D_X(M)\cong f_*\RHom_{\Omega_X}(M,\Omega_X[2\dim_S X])&\to  \\
\RHom_{\Omega_Y}(f_*M,f_*\Omega_X[2\dim_S X])&\xrightarrow{\Tr_p} \RHom_{\Omega_Y}(f_*M,\Omega_Y[2\dim_S Y])\cong D_Y(f_*M).
\end{split}\end{equation*}
It is enough to prove that this morphism is an isomorphism when $Y$ is affine, and $M\cong q^*\OO_{\PP^n}(k)$, where $q:X\to \PP^n$ is the second projection. This is straightforward.
\qed

\begin{thm}\label{adj_prop_smooth}
Suppose that $f:X\to Y$ is smooth and proper. Let $N\in\D^\co(\qc{X})$, and let $M\in\D^\abs(\coh{X})$. Then
$$\RHom_{\Omega_Y}(f_*M,N))\cong f_*\RHom_{\Omega_X}(M,f^!N).$$ 
\end{thm}
\proof
This follows from the projection formula, Proposition~\ref{dual_hom}, and Theorem~\ref{dual_direct}:
\begin{equation*}\begin{split}
&f_*\RHom_{\Omega_X}(M,f^!N)\cong f_*(D(M)\otimes^D_{\Omega_X}f^!N)\cong f_*D(M)\otimes^D_{\Omega_Y}N\cong \\
&D(f_*M)\otimes^D_{\Omega_Y}N\cong \RHom_{\Omega_Y}(f_*M,N)
\end{split}\end{equation*}
\qed

The proof of the following result is the same as the proof of~\cite[2.7.1]{HTT}

\begin{prop}\label{dual_inverse}
If $f:X\to Y$ is smooth, then $D_X\circ f^!\cong f^!\circ D_Y[ 2(\dim X-\dim Y)]$. Equivalently, 
$f^+\cong D_X\circ f^!\circ D_Y$.
\end{prop}

\section{$\Omega$-modules and D-modules.}\label{comparison}

\subsection{Positselski theorem.} Suppose that $S=\Spec k$, where $k$ is a filed of characteristic zero. 
Let $\DD_X$ be the sheaf of differential operators on a smooth quasi-projective algebraic variety $X$ over $k$. We have an $\Omega_X-\DD_X$--bimodule
$\cDR_X=\Omega_X\otimes_{\OO_X}\DD_X$ with differential given in local coordinates $x_1,\dots,x_n$ on $X$ by the formula:
$$d(\omega\otimes m)=d\omega\otimes m+\sum_{i=1}^n \omega\wedge dx_i\otimes \d/\d x_im,$$ where $m$ is a local section of $\DD_X$ and $\omega$ is a local section of $\Omega_X$. The following result is an easy corollary of~\cite[B.2,B.5]{P2}. 

\begin{thm*}\cite[B.2,B.5]{P2}
The functor $\DR_X=-\otimes_{\DD_X}\cDR_X[\dim X]:\Dmod{X}\to\qc{X}$ induce an equivalence of the unbounded derived category of quasi-coherent left $\DD_X$-modules $\D(\Dmod{X})$ and the coderived category $\D^\co(\qc{X})$. Moreover, this functor induce an equivalence of the bounded derived category of coherent left $\DD_X$-modules $\D^b(\Dcoh{X})$ and the absolute derived category $\D^\abs(\coh{X})$.
\end{thm*}

\begin{lemma}\label{coiso}
Let $\cDR_X^{-1}=\DD_X\otimes_{\OO_X}\omega_X^{-1}\otimes_{\OO_X}\Omega_X$ be a $\DD_X-\Omega_X$--bimodule with differential $$d(m\otimes\theta_1\wedge\dots\wedge\theta_r)=\sum_{i=1}^rm\theta_i\otimes \theta_1\wedge\dots\wedge\Hat\theta_i\wedge\dots\wedge\theta_r,$$
where $\theta_1,\dots,\theta_r$ are commuting vector fields.
%$\mathrm{(}$see~\cite[1.5.27]{HTT} for the definition of differential$\mathrm{)}$
Then there exists a natural co-isomorphism of $\Omega_X$-modules $\OO_X\to\cDR_X^{-1}[\dim X].$
\end{lemma}
\proof
The morphism is induced by the natural inclusion $\OO_X\to \DD_X$. The filtration of $\cDR_X^{-1}/\OO_X$ defined in~\cite[1.5.27]{HTT} has acyclic associated graded complex. Each graded component is a bounded complex of $\OO_X$-modules and thus coacyclic. The Lemma follows from Lemma~\ref{colemma}.(3).
\qed

\subsection{} 
Let $f:X\to Y$ be a morphism of smooth quasi-projective algebraic varieties over $k$. Recall the definition of the direct image $f_\star$ and the extraordinary inverse image $f^\dagger$ for $\DD$-modules~\cite{HTT}. 
The functor $f^\dagger$ is the derived functor of $f^\bullet(-)\otimes_{f^\bullet\OO_Y}\OO_X$ shifted by $\dim X-\dim Y$. Let $\DD_{Y\leftarrow X}=\omega_X\otimes_{f^\bullet\OO_Y}f^\bullet(\DD_Y\otimes_{\OO_Y}\omega_Y^{-1})$ be a $\DD_X-f^\bullet \DD_Y$--bimodule. Then $f_\star=Rf_\bullet(\DD_{Y\leftarrow X}\otimes^L_{\DD_X}-)$. For coherent $\DD_X$-modules we have the duality functor $\D'_X(-)=\RHom_{\DD_X}(-,\DD_X\otimes_{\OO_X}\omega_X^{-1})[\dim X]$.
We are going to prove that the Positselski equivalence takes these functors to corresponding functors for $\Omega$-modules. We also prove that the Positselski equivalence is a tensor functor. By definition, the tensor product of left $\DD_X$-modules $M$ and $N$ is the left $\DD_X$-module $M\otimes_{\OO_X}N$ such that a vector field $\theta$ acts by the formula: $\theta(m\otimes n)=\theta(m)\otimes n+m\otimes\theta(n)$.

\begin{thm*}\label{comp}
\begin{enumerate}
\item There exists an equivalence of bifunctors from 
$\D(\Dmod{X})\times\D(\Dmod{X})$ to $\D^\co(\qc{X})$: $$\DR_X(-)\otimes^D_{\Omega_X}\DR_X(-)\cong\DR_X(-\otimes^L_{\OO_X}-).$$
\item $f_*(\DR_X(-))\cong\DR_Y(f_\star(-))$.
\item $f^!(\DR_Y(-))\cong\DR_X(f^\dagger(-))$.
\item $\D_X(\DR_X(-))\cong\DR_X(\D'_X(-)).$
\end{enumerate}
\end{thm*}
\proof
\begin{enumerate}
\item Let $M$ and $N$ be complexes of quasi-coherent $\DD_X$-modules. Assume that $M$ is a complex of flat modules over $\OO_X$. We have the following co-isomorphisms in $\D^\co(\qc{X})$:
$$\DR_X(M)\otimes^D_{\Omega_X}\DR_X(N)\cong (\Omega_X\otimes_{\OO_X}M)\otimes_{\Omega_X}(\Omega_X\otimes_{\OO_X}N)\cong
\Omega_X\otimes_{\OO_X}(M\otimes_{\OO_X}N)\cong \DR_X(M\otimes^D_{\OO_X}N).$$

\item Let $L$ be a complex of quasi-coherent flat $\DD_Y$-modules. Then
\begin{equation*}\begin{split}
&f^!(\DR_Y(L))\cong \Omega_X\otimes_{f^\bullet\Omega_Y}f^\bullet(\Omega_Y\otimes_{\OO_Y}L)[2\dim X-\dim Y]\cong \\
&\Omega_X\otimes_{f^\bullet\OO_Y}f^\bullet L[2\dim X-\dim Y]\cong 
\Omega_X\otimes_{\OO_X}f^\dagger L[\dim X]\cong \DR_X(f^\dagger(L)).
\end{split}\end{equation*}

\item We give two proofs of this fact. First, the morphism $f$ is a composition of a closed embedding $i:X\to Z$ and a smooth projection $g:Z\to Y$. For example, $Z=X\times_S Y$, and $i$ is the graph of $f$. Clearly, $(3)$ follows from $(2)$ and Theorems~\ref{inj_closed} and~\ref{adj_smooth}, since the adjoint functor is unique up to equivalence. Unfortunately, this equivalence depends on the choise of $g$ and $i$, and we give a more explicit second proof. 

Let $M$ be a complex of flat quasi-coherent $\DD_X$-modules. Then 
\begin{equation*}\begin{split}
&\DR_Y(f_\star(M))\cong\Omega_Y\otimes_{\OO_Y}Rf_\bullet(\DD_{Y\leftarrow X}\otimes_{\DD_X}M)[\dim Y]\cong \\
&\Omega_Y\otimes_{\OO_Y}Rf_\bullet(f^\bullet(\omega_Y^{-1}\otimes_{\OO_Y}\DD_Y)\otimes_{f^\bullet\OO_Y} \omega_X\otimes_{\DD_X}M)[\dim Y]\cong \\
&\Omega_Y\otimes_{\OO_Y}Rf_\bullet(f^\bullet(\omega_Y^{-1}\otimes_{\OO_Y}\DD_Y)\otimes_{f^\bullet\OO_Y}\Omega_X\otimes_{\OO_X}\DD_X\otimes_{\DD_X}M)[\dim X+\dim Y]\cong \\
&\Omega_Y\otimes_{\OO_Y}Rf_\bullet(f^\bullet(\omega_Y^{-1}\otimes_{\OO_Y}\DD_Y)\otimes_{f^\bullet\OO_Y}\Omega_X\otimes_{\OO_X}M)[\dim X+\dim Y],
\end{split}\end{equation*}
and
\begin{equation*}\begin{split}
&f_*(\DR_X(M))\cong \Omega_Y\otimes_{\Omega_Y}f_*(\Omega_X\otimes_{\OO_X}M)[\dim X]\cong\\ 
&\Omega_Y\otimes_{\OO_Y}\omega^{-1}_Y\otimes_{\OO_Y}\DD_Y\otimes_{\OO_Y}\Omega_Y\otimes_{\Omega_Y}f_*(\Omega_X\otimes_{\OO_X}M)[\dim X+\dim Y],
\end{split}\end{equation*}
where we used the quasi-isomorphism of right $\DD_X$-modules $\Omega_X\otimes_{\OO_X}\DD_X[\dim X]\to\omega_X$ (see~\cite[1.5.27]{HTT}), and the corresponding co-isomorphism of $\Omega_Y$-modules $\Omega_Y\otimes_{\OO_Y}\omega^{-1}_Y\otimes_{\OO_Y}\DD_Y\otimes_{\OO_Y}\Omega_Y\to\Omega_Y[\dim Y]$.

Choose an affine covering $\UU$ of $X$. We can compute $f_*$ and $Rf_\bullet$ using corresponding \v{C}ech resolutions. Thus it is enough to construct the equivalence when $X$ and $Y$ are affine. We have to check that the natural morphism 
\begin{equation*}\begin{split}
&\Omega_Y\otimes_{\OO_Y}\omega^{-1}_Y\otimes_{\OO_Y}\DD_Y\otimes_{\OO_Y}\Omega_Y\otimes_{\Omega_Y}\Omega_X\otimes_{\OO_X}M\to\\ 
&\Omega_Y\otimes_{\OO_Y}\omega_Y^{-1}\otimes_{\OO_Y}\DD_Y\otimes_{\OO_Y}\Omega_X\otimes_{\OO_X}M
\end{split}\end{equation*}
 commutes with differentials. This follows from definitions.

\item
The functor $\DR_X^{-1}(-)=\cDR_X^{-1}\otimes_{\Omega_X}-$ is quasi-inverse to $\DR_X$.
Let $P$ be a perfect $\Omega_X$-module. We prove that $\D'_X(\DR^{-1}(P))\cong\DR_X^{-1}(\D_X(P))$. 
Using the quasi-isomorphism of left $\DD_X$-modules  
$\DD_X\otimes_{\OO_X}\omega_X^{-1}\otimes_{\OO_X}\Omega_X[\dim X]\to \OO_X$ we get:
 
\begin{equation*}\begin{split} 
&\D'_X(\DR^{-1}(P))\cong \cHom_{\DD_X}(\DD_X\otimes_{\OO_X}\omega_X^{-1}\otimes_{\OO_X}P,\DD_X\otimes_{\OO_X}\omega_X^{-1})[\dim X]\cong \\ 
&\cHom_{\DD_X}(\DD_X\otimes_{\OO_X}\omega_X^{-1}\otimes_{\OO_X}P,\DD_X\otimes_{\OO_X}\omega_X^{-1} \otimes_{\OO_X}\OO_X[-\dim X])[2\dim X]\cong \\ 
&\cHom_{\DD_X}(\DD_X\otimes_{\OO_X}\omega_X^{-1}\otimes_{\OO_X}P,\DD_X\otimes_{\OO_X}\omega_X^{-1}\otimes_{\OO_X}\DD_X\otimes_{\OO_X}\omega_X^{-1}\otimes_{\OO_X}\Omega_X)[2\dim X],
\end{split}\end{equation*}
We have the natural transformation
\begin{equation*}\begin{split}
&\Phi:\DR_X^{-1}(\D_X(P))[-2\dim X]\cong \DD_X\otimes_{\OO_X}\omega_X^{-1}\otimes_{\OO_X}\cHom_{\Omega_X}(P,\Omega_X)\to\\
&\cHom_{\DD_X}(\DD_X\otimes_{\OO_X}\omega_X^{-1}\otimes_{\OO_X}P,\DD_X\otimes_{\OO_X}\omega_X^{-1}\otimes_{\OO_X}\DD_X\otimes_{\OO_X}\omega_X^{-1}\otimes_{\OO_X}\Omega_X)
\end{split}\end{equation*}
given over the affine subscheme $U\subset X$ by the formula $m\otimes \phi\mapsto(d\otimes p\mapsto(m\otimes d\otimes\phi(p))$, where $\phi\in\cHom_{\Omega_X}(P,\Omega_X)(U)$, $p\in P(U)$, $m,d\in \DD_U\otimes_{\OO_U}\omega_U^{-1}$. We claim that the transformation is an isomorphism of $\DD_X$-modules. As usual, it is enough to prove the claim locally and for $P=\OO_X$. This is clear.
\end{enumerate}
\qed

\begin{rem*}
The internal Hom for $\Omega_X$-modules corresponds to the internal Hom for $\DD$-modules under Positselski equivalence. The last bifunctor is defined as follows: $$\RHom_{\DD_X}(M,N)=\D_X'(M)\otimes^L_{\OO_X}N,$$ where $M$ is coherent, and $N$ is quasi-coherent. 
\end{rem*}

%\section{Open questions.}
%\subsesction{}
%\subsection{Constractible $\Omega_X$-modules.}

\end{document}